\documentclass[reqno,english]{amsart}
\usepackage[T1]{fontenc}
\usepackage{lmodern,mathtools,amssymb,microtype}
\usepackage[colorlinks=true,hypertexnames=false,
 linkcolor=blue,citecolor=blue,urlcolor=blue]{hyperref}

\numberwithin{equation}{section}
\newtheorem{theorem}{Theorem}[section]
\newtheorem{lemma}{Lemma}[section]
\newtheorem{proposition}{Proposition}[section]
\newtheorem{corollary}{Corollary}[section]

\let\bccbibliography\thebibliography
\renewcommand{\thebibliography}[1]{%
 \bccbibliography{#1}%
 \setlength{\itemsep}{0pt}%
}

\newcommand{\R}{\mathbb R}
\newcommand{\Z}{\mathbb Z}
\newcommand{\E}{\mathcal E}
\newcommand{\tr}{\operatorname{tr}}
\newcommand{\dist}{\operatorname{dist}}
\newcommand{\diag}{\operatorname{diag}}
\newcommand{\Span}{\operatorname{span}}
\newcommand{\SO}{\operatorname{SO}}
\allowdisplaybreaks[2]

\title[BCC quantizer stability]
 {\bfseries On sharp stability of the BCC quantizer}
\author[S.~Luo]{Senping Luo}
\author[B.~Xie]{Bingbing Xie}

\address[S.~Luo and B.~Xie]{School of Mathematics and statistics, Jiangxi Normal University, Nanchang, 330022, China; and Jiangxi applied mathematical research center.}

\email[S.~Luo]{luosp1989@163.com }
\email[B.~Xie]{xiebb2074@163.com }

\date{}
\hypersetup{
 pdftitle={On sharp stability of the BCC quantizer},
 pdfauthor={Senping Luo and Bingbing Xie}
}

\begin{document}

\begin{abstract}
We prove sharp nonlinear stability estimates for the squared-distance quantization energy of the body-centered cubic lattice. The reference Hessian has optimal coercivity constant \(1/96\) relative to nearest-neighbor displacement differences, uniformly over all finite-index periods. We classify all equality modes and prove cubic cancellation on their span. This yields nonlinear stability in a neighborhood independent of the period, including volume-preserving strain. The loss in the sharp energy coefficient is at most quadratic for each fixed period and exactly linear when all periods are allowed. We also determine the transition between these regimes. Finally, we prove global optimality and uniqueness of BCC, up to rigid motions, among two-site configurations with rectangular period cells of fixed volume.
\end{abstract}

\maketitle

\section{Introduction and main results}
The squared-distance quantization problem asks how a set of generators should be arranged to minimize the distance to the nearest generator. We consider periodic configurations in dimension three and perturbations of the body-centered cubic lattice (BCC).

Let
\begin{equation}\label{eq:lattice}
 \Lambda=\Z^3\cup\bigl(\Z^3+(\tfrac12,\tfrac12,\tfrac12)\bigr),
 \qquad
 \mathcal R=\{\tfrac12(1,\sigma,\tau):\sigma,\tau\in\{-1,1\}\}.
\end{equation}
The covolume of \(\Lambda\) is \(1/2\). Let \(\Gamma\subset\Lambda\) be a finite-index sublattice, let \(N=[\Lambda:\Gamma]\), and let \(\Omega_\Gamma\) be a fundamental domain of \(\Gamma\). For a \(\Gamma\)-periodic field \(u:\Lambda\to\R^3\), put
\begin{align}
 \Lambda_u&=\{p+u_p:p\in\Lambda\},&
 E_\Gamma(u)&=\int_{\Omega_\Gamma}\dist(x,\Lambda_u)^2\,dx,\label{eq:energy}\\
 D_\Gamma(u)&=\sum_{p\in\Lambda/\Gamma}\sum_{r\in\mathcal R}
 |u_{p+r}-u_p|^2,&
 \eta(u)&=\max_{\substack{p\in\Lambda/\Gamma\\r\in\mathcal R}}
 |u_{p+r}-u_p|.\label{eq:gradient}
\end{align}
All edges retain their translation labels. On a small quotient, distinct labeled edges may have the same endpoints or may become loops. The definitions retain their multiplicities.

The Voronoi cell of \(\Lambda\) is a truncated octahedron of volume \(1/2\), with squared-distance moment \(19/256\). Hence \(E_\Gamma(0)=19N/256\). The integral in \eqref{eq:energy} is the classical quantization energy. At configurations with distinct generators, its stationary points are centroidal Voronoi tessellations; see Du, Faber and Gunzburger \cite{DuFaberGunzburger}. General quantization theory and its geometric applications are developed in \cite{GrafLuschgy,Gruber}.

Barnes and Sloane \cite{BarnesSloane} proved that BCC is the optimal lattice quantizer in dimension three; see also Conway and Sloane \cite{ConwaySloane}. Optimality among arbitrary point configurations is a stronger problem, related to Gersho's conjecture \cite{Gersho}, and remains open \cite[Section 1.3]{BourneSchmitzerWirth}. Bounds on the geometric complexity of optimal three-dimensional cells were obtained by Choksi and Lu \cite{ChoksiLu}. Numerical evidence for BCC is provided in \cite{DuWang,Hain,Klatt}. In two dimensions, Tisdell, Choksi and Lu \cite{TisdellChoksiLu} recently established quantitative bounds on the number of nonhexagonal cells for optimal quantizers on the hexagonal torus. Quantization and its gradient-flow formulation are studied in \cite{Caglioti}.

The first and second variations of Voronoi energies have been studied in \cite{Urschel}. For the BCC energy, Ruscher, Baschnagel and Farago \cite[Section III.F and Appendix B]{Ruscher} give a Hessian representation with longitudinal and transverse face contributions. Periodic centroidal tessellations are considered in \cite{Zhang}, and numerical stability under small BCC perturbations appears for related location energies in \cite[online supplement, Section 1.2]{BourneRoper}. A general Fourier framework for lattice stability and nonlinear perturbations is developed in \cite{HudsonOrtner,OrtnerTheil}. These works provide the variational setting for the sharp estimates below. Our analysis determines the exact coercivity coefficient, all equality modes, and the joint dependence of the nonlinear correction on the period and the perturbation size.

Quantization also has a ground-state interpretation. Torquato \cite{Torquato} relates the quantizer problem to particle interactions of several orders. The same many-body energy is studied as the Voronoi liquid in \cite{Ruscher,RuscherSound,Hain}.

Another mechanism selecting BCC appears in the competing interaction model of Ren and Wei \cite{RenWeiBCC}, motivated by diblock copolymers and combining perimeter with a nonlocal interaction generated by \((-\Delta)^{-3/2}\). They proved that the regular part of the corresponding periodic Green function at the origin is uniquely minimized by BCC among three-dimensional lattices of unit covolume, using its relation to the height of the dual flat torus.

More generally, lattice selection is a central question in crystallization; see Blanc and Lewin \cite{BlancLewin}. For pair-interaction energies, the first author and Wei proved sharp minimization results for differences of theta functions \cite{LuoWeiTheta} and Epstein zeta functions, including Lennard--Jones lattice energies \cite{LuoWeiZeta}. Theta derivatives and lattice transitions are treated in \cite{LuoWeiPhase}, and their recent three-dimensional work \cite{LuoWeiSS} concerns theta and Epstein zeta optimization among Bravais lattices. Here the energy is determined by the nearest-site partition, and the admissible perturbations include nonlattice configurations.

The main purpose of this paper is to establish sharp stability of BCC uniformly over all finite-index periods and to determine the associated nonlinear correction. We also prove global optimality in the rectangular two-site class. Our first result is the following.

\begin{theorem}[Sharp coercivity and uniform nonlinear stability]\label{thm:stability}
For every finite-index period sublattice \(\Gamma\) and every periodic field \(v\),
\begin{equation}\label{eq:sharp-main}
 D^2E_\Gamma(0)[v,v]\geq\frac1{96}D_\Gamma(v).
\end{equation}
The constant \(1/96\) is optimal uniformly over all periods; every equality frequency and polarization is classified in Proposition~\ref{prop:equality}. There exists \(\eta_0>0\), independent of \(\Gamma\), such that
\begin{equation}\label{eq:nonlinear-main}
 \eta(u)\leq\eta_0
 \quad\Longrightarrow\quad
 \begin{cases}
 D^2E_\Gamma(u)[v,v]\geq D_\Gamma(v)/192,\\
 E_\Gamma(u)-E_\Gamma(0)\geq D_\Gamma(u)/384.
 \end{cases}
\end{equation}
Here \(E_\Gamma(0)=19N/256\). In this neighborhood, \(E_\Gamma(u)=E_\Gamma(0)\) if and only if \(u\) is constant, and every periodic stationary configuration is a translate of BCC.
\end{theorem}

The radius in Theorem~\ref{thm:stability} is independent of both the period and the number of sites. It controls displacement differences, and therefore permits displacements whose amplitude grows with the period. The quadratic Taylor term associated with \eqref{eq:sharp-main} has energy coefficient \(1/192\). Our next result determines the order of the nonlinear correction to this coefficient.

To formulate the loss precisely, set
\begin{equation}\label{eq:coefficient-defs}
 c_\Gamma(\rho)=\inf_{0<\eta(u)\leq\rho}
 \frac{E_\Gamma(u)-E_\Gamma(0)}{D_\Gamma(u)},
 \qquad
 c_{\mathrm{all}}(\rho)=\inf_{\Gamma\subset\Lambda}c_\Gamma(\rho).
\end{equation}
The second infimum runs over finite-index sublattices. Quotients admitting only translations are omitted from these infima.

\begin{theorem}[Fixed-period and all-period loss]\label{thm:loss}
For every fixed \(\Gamma\), there is a finite constant \(C_\Gamma\) such that, for all sufficiently small \(\eta(u)\),
\begin{equation}\label{eq:fixed-main}
 E_\Gamma(u)-E_\Gamma(0)\geq
 \left(\frac1{192}-C_\Gamma\eta(u)^2\right)D_\Gamma(u).
\end{equation}
For \(\Gamma=4\Z^3\), the loss has exactly quadratic order:
\begin{equation}\label{eq:four-main}
 \frac{\rho^2}{15744}\leq
 \frac1{192}-c_{4\Z^3}(\rho)\leq C_{4\Z^3}\rho^2.
\end{equation}
Across all periods, it has exactly linear order. More precisely, for all sufficiently small \(\rho>0\),
\begin{equation}\label{eq:all-main}
 \frac{\rho}{400\sqrt2}\leq
 \frac1{192}-c_{\mathrm{all}}(\rho)
 \leq\frac{10\sqrt2M_3}{3}\rho,
\end{equation}
where \(M_3\) is the universal local derivative bound in Lemma~\ref{lem:local-cell}.
\end{theorem}

The difference between the two orders in Theorem~\ref{thm:loss} is determined by the equality modes of the reference Hessian. The cubic variation vanishes on their full span. On a fixed quotient, the remaining frequencies are separated from this space, giving the quadratic correction in \eqref{eq:fixed-main}. This separation decreases as the period grows. We determine its precise effect in the following theorem.

For \(\Gamma_m=m\Z^3\), write \(E_m=E_{\Gamma_m}\), \(D_m=D_{\Gamma_m}\), and \(c_m=c_{\Gamma_m}\). Let \(C_m^{\mathrm{loc}}\) be the infimum of all constants \(A\) for which
\[
 E_m(u)-E_m(0)\geq
 \left(\frac1{192}-A\eta(u)^2\right)D_m(u)
\]
holds on some positive neighborhood \(\eta(u)\leq\rho_m\).

\begin{theorem}[Sharp period dependence and crossover]\label{thm:crossover}
There are constants \(C,\rho_*>0\), independent of \(m\), such that
\begin{equation}\label{eq:period-main}
 E_m(u)-E_m(0)\geq
 \left(\frac1{192}-Cm^2\eta(u)^2\right)D_m(u),
 \qquad \eta(u)\leq\rho_*.
\end{equation}
Moreover,
\begin{equation}\label{eq:growth-main}
 \liminf_{\substack{m\to\infty\\4\mid m}}
 \frac{C_m^{\mathrm{loc}}}{m^2}\geq\frac1{342\pi^2},
 \qquad C_m^{\mathrm{loc}}\leq Cm^2.
\end{equation}
There exist universal \(c,C,\rho_0>0\) and \(m_0\geq8\) such that
\begin{equation}\label{eq:crossover-main}
 c\min\{\rho,m^2\rho^2\}\leq
 \frac1{192}-c_m(\rho)
 \leq C\min\{\rho,m^2\rho^2\}
\end{equation}
whenever \(4\mid m\), \(m\geq m_0\), and \(0<\rho\leq\rho_0\).
Thus the transition occurs at \(\rho\asymp m^{-2}\).
\end{theorem}

We next allow the period lattice itself to vary. Write
\[
 F=Re^S,\qquad R\in\SO(3),\qquad S=S^T,\qquad \tr S=0,
\]
and define
\begin{equation}\label{eq:strained-energy}
 \Lambda_{F,u}=\{Fp+u_p:p\in\Lambda\},
 \qquad
 \E_\Gamma(F,u)=\int_{F\Omega_\Gamma}
 \dist(x,\Lambda_{F,u})^2\,dx.
\end{equation}
This configuration has period lattice \(F\Gamma\) and unchanged density.

\begin{theorem}[Stability under volume-preserving strain]\label{thm:strain}
There is a universal \(\rho>0\) such that
\begin{equation}\label{eq:strain-main}
 \E_\Gamma(F,u)-\frac{19N}{256}
 \geq\frac{N}{192}|S|_F^2+\frac1{384}D_\Gamma(u),
 \qquad |S|_F+\eta(u)\leq\rho.
\end{equation}
Equality with the BCC energy holds exactly when \(S=0\) and \(u\) is constant. The homogeneous strain expansion is
\begin{equation}\label{eq:strain-expansion}
 \frac{\E_\Gamma(e^S,0)}N
 =\frac{19}{256}+\frac1{96}\sum_iS_{ii}^2
 +\frac1{16}\sum_{i<j}S_{ij}^2+O(|S|_F^3).
\end{equation}
The remainder constant is independent of the period.
\end{theorem}

Theorems~\ref{thm:stability}--\ref{thm:strain} allow every site in the period cell to move independently. The assumptions concern only the displacement differences and the homogeneous strain. We now consider a different question: global minimization over all configurations with two sites per rectangular period cell. In this class, the relative shift and the aspect ratios can be arbitrary.

\begin{theorem}[Global optimality for rectangular two-site periods]\label{thm:global}
Let \(\ell_i>0\), \(\ell_1\ell_2\ell_3=1\), and \(a\in\R^3/\Z^3\). Set
\[
 D_\ell=\diag(\ell_1,\ell_2,\ell_3),\qquad
 L_\ell(a)=D_\ell\bigl(\Z^3\cup(\Z^3+a)\bigr),\qquad
 \Omega_\ell=D_\ell[0,1]^3.
\]
Then
\begin{equation}\label{eq:global-main}
 Q(\ell,a):=\int_{\Omega_\ell}\dist(x,L_\ell(a))^2\,dx
 \geq\frac{19}{128}.
\end{equation}
Equality holds if and only if \(\ell_1=\ell_2=\ell_3=1\) and
\(a=(1/2,1/2,1/2)\pmod{\Z^3}\). Thus, up to a rigid motion, BCC is the unique global minimizer in this class. For a cubic period, writing
\(\delta_i=\dist(a_i,1/2+\Z)\), one has the sharp bound
\begin{equation}\label{eq:global-gap-main}
 Q(\mathbf1,a)-\frac{19}{128}
 \geq\frac7{96}\sum_{i=1}^3\delta_i^2,
 \qquad \mathbf1=(1,1,1).
\end{equation}
The coefficient \(7/96\) is the largest possible.
\end{theorem}

Theorem~\ref{thm:global} permits changes of Voronoi incidences and includes point sets that are not Bravais lattices. Proposition~\ref{prop:global-rectangular} gives a quantitative estimate for the rectangular deformation, and Proposition~\ref{prop:global-cubic} strengthens \eqref{eq:global-gap-main} to an explicit polynomial inequality. The competitor class is restricted to rectangular periods with two sites. The unrestricted BCC quantization problem, with arbitrary periods and arbitrarily many sites per period, remains open.

We conclude the introduction with the main reductions in the proof. For sufficiently small \(\eta(u)\), the energy has the local representation
\[
 E_\Gamma(u)=\sum_{p\in\Lambda/\Gamma}\Phi(d_p(u)),
 \qquad d_p(u)=(u_{p+r}-u_p)_{r\in\mathcal N},
\]
where \(\mathcal N=\mathcal R\cup(-\mathcal R)\cup\{\pm e_1,\pm e_2,\pm e_3\}\). The analytic function \(\Phi\) is determined by these fourteen neighbors. The same function applies to every period. The Fourier symbol \(H(k)\) of the reference Hessian satisfies
\[
 H(k)=\frac1{96}M(k/2),\qquad
 M(x)\succeq8(1-\cos x_1\cos x_2\cos x_3)I.
\]
The matrix \(M\) is given explicitly in \eqref{eq:symbol}. Two polynomial identities prove this inequality and identify every equality frequency and polarization. This establishes the sharp reference estimate and, through the local derivative bounds, the nonlinear stability in Theorem~\ref{thm:stability}.

For the nonlinear correction, let \(K_\Gamma\) be the equality space after removing translations and put
\[
 g_\Gamma(u)=D^2E_\Gamma(0)[u,u]-\frac1{96}D_\Gamma(u).
\]
The decisive cancellation is
\[
 D^3E_\Gamma(0)[v_1,v_2,v_3]=0
 \qquad(v_1,v_2,v_3\in K_\Gamma).
\]
It follows from the reciprocal-lattice resonance relations and the signed coordinate symmetries of BCC. With orthogonality taken in the \(D_\Gamma\) inner product, write \(g_m=g_{\Gamma_m}\) and \(K_m=K_{\Gamma_m}\). For a constant \(b>0\) independent of \(m\), the transverse estimate
\[
 g_m(w)\geq\frac{b}{m^2}D_m(w),\qquad w\in K_m^\perp,
\]
and a uniform projection bound give the upper estimate in \eqref{eq:crossover-main}. An explicit coupling of an equality mode to a nearby frequency gives the matching lower estimate. The balance between the quadratic transverse cost and the cubic coupling occurs at \(\rho\asymp m^{-2}\).

The global argument uses a separate exact reduction:
\begin{align*}
 Q(\ell,a)&=\frac{\ell_1^2+\ell_2^2+\ell_3^2}{12}-\frac12A(w),
 \qquad w_i=\ell_i^2a_i(1-a_i),\\
 A(w)&=\frac18\int_{[-1,1]^3}|w_1v_1+w_2v_2+w_3v_3|\,dv,
\end{align*}
where \(a_i\in[0,1]\). A lower bound for each width derivative of \(A\) first reduces the minimization to the centered shift. The remaining aspect-ratio problem reduces to the strictly decreasing function \(f\) in \eqref{eq:global-shape-profile}. Thus the global minimum is attained only at the cubic shape.

Section~\ref{sec:geometry} establishes the local cell representation and the exact Hessian. Section~\ref{sec:fourier} proves sharp coercivity, classifies the equality modes, and establishes cubic cancellation, uniform nonlinear stability, and the strain estimate. Section~\ref{sec:fixed} determines the nonlinear loss, from the fixed-period expansion and four-phase construction to the sharp period dependence and crossover law. Section~\ref{sec:global} proves global optimality for rectangular two-site periods and classifies the stationary shifts for a cubic period. The appendices contain the polynomial coefficient estimates and exact local cell identities.

\section{Voronoi geometry and the exact second variation}\label{sec:geometry}
In this section, we compute the second variation of the quantization energy at BCC. The following exact bond formula is the main result.

\begin{proposition}[Exact bond representation]\label{prop:bond}
The first variation vanishes at BCC, and its Hessian is
\begin{equation}\label{eq:bond}
 D^2E_\Gamma(0)[u,u]=
 \sum_{p\in\Lambda/\Gamma}\sum_{r\in\mathcal R\cup\{e_1,e_2,e_3\}}
 (u_{p+r}-u_p)^TB_r(u_{p+r}-u_p),
\end{equation}
where
\begin{equation}\label{eq:bond-matrices}
 B_{e_i}=\frac{25e_ie_i^T-I}{384},
 \qquad B_{\sigma/2}=\frac{41\sigma\sigma^T-15I}{768},
 \quad \sigma\in\{1\}\times\{-1,1\}^2.
\end{equation}
\end{proposition}

The proof is based on the geometry of the reference Voronoi cell and the variation of its interfaces. We first establish a local analytic representation of the energy and compute the required face moments. All estimates involve only local displacement differences, which makes their constants independent of the period.

\begin{lemma}[Local analytic cell representation]\label{lem:local-cell}
There exist \(\delta>0\) and \(M_3,M_4\geq1\), depending only on the reference BCC cell, with the following property. Let
\[
 \mathcal N=\mathcal R\cup(-\mathcal R)\cup\{\pm e_1,\pm e_2,\pm e_3\},
 \qquad
 K(d)=\bigcap_{r\in\mathcal N}
 \{y:(r+d_r)\cdot y\leq |r+d_r|^2/2\}.
\]
For \(\|d\|\leq\delta\), this polytope has the reference BCC face incidences, and
\[
 \Phi(d)=\int_{K(d)}|y|^2\,dy
\]
is analytic on a neighborhood of the closed parameter ball, with
\(\|D^j\Phi(d)\|_{\mathrm{op}}\leq M_j\), \(j=3,4\). If
\[
 \eta(u)\leq\min\{1/100,\delta/(4\sqrt2)\},
 \qquad d_p(u)=(u_{p+r}-u_p)_{r\in\mathcal N},
\]
then the actual Voronoi cell at \(p+u_p\), translated to its generator, is \(K(d_p(u))\), and
\begin{equation}\label{eq:local-energy}
 E_\Gamma(u)=\sum_{p\in\Lambda/\Gamma}\Phi(d_p(u)).
\end{equation}
\end{lemma}

We recall that the reference Voronoi cell is
\begin{equation}\label{eq:reference-cell}
 V=\{y:|y_i|\leq\tfrac12,\quad |y_1|+|y_2|+|y_3|\leq\tfrac34\}.
\end{equation}
Its vertices are the permutations of \((0,\pm1/4,\pm1/2)\), its volume is \(1/2\), its covering radius is \(R_0=\sqrt5/4\), and
\begin{equation}\label{eq:cell-moment}
 \int_V|y|^2\,dy=\frac{19}{256},
 \qquad \int_Vyy^T\,dy=\frac{19}{768}I.
\end{equation}
In the standard normalization, a cell of volume \(v\) has dimensionless quantization constant \(G=(3v^{5/3})^{-1}\int_V|y|^2\,dy\). Thus \eqref{eq:cell-moment} agrees with the BCC value in \cite{BarnesSloane,ConwaySloane}. We shall also use the following face moments.

\begin{lemma}[Exact face moments]\label{lem:face-moments}
Let $F_r$ be the reference face shared with the neighbor $r=dn$, where
$|n|=1$, and write $y=x-r/2$ on this face. Then
\[
 \int_{F_r}yy^T\,dS=A\beta(I-nn^T),\qquad
 (A,d,\beta)=
 \begin{cases}
 (1/8,1,1/96),&r=\pm e_i,\\
 (3\sqrt3/16,\sqrt3/2,5/192),
   &r=\tfrac12(\pm1,\pm1,\pm1).
 \end{cases}
\]
\end{lemma}
\begin{proof}
The two types of faces are a square and a regular hexagon, both with
side length $a=\sqrt2/4$. By rotational symmetry, the centered
second moment is scalar on the tangent plane. For the square, integration
in either tangential coordinate gives $\beta=a^2/12$. For the hexagon,
we divide it into six triangles with a common vertex at the center. Each triangle has
two radius vectors of length $a$ and mutual scalar product $a^2/2$;
integration of the barycentric quadratic monomials therefore gives
$\int |y|^2\,dS=5a^2/12$ times the triangle's area. Summing over the
triangles and dividing by the two tangent dimensions yields
$\beta=5a^2/24$. The polygon areas and neighbor lengths give the stated
constants.
\end{proof}

To identify the perturbed cells using only these fourteen neighbors, we first estimate the displacement difference between arbitrary lattice sites.

\begin{lemma}[Nearest-neighbor graph distance]\label{lem:distance}
The word length of \(v\in\Lambda\) in the eight steps \(\mathcal R\cup(-\mathcal R)\) is \(2\|v\|_\infty\). Consequently
\begin{equation}\label{eq:displacement-distance}
 |u_q-u_p|\leq2\eta(u)\|q-p\|_\infty
 \leq2\eta(u)|q-p|.
\end{equation}
In particular \(p\mapsto p+u_p\) is injective for \(\eta(u)<1/2\), and
\begin{equation}\label{eq:bilipschitz}
 (1-2\eta)|q-p|\leq |q+u_q-p-u_p|
 \leq(1+2\eta)|q-p|.
\end{equation}
\end{lemma}

\begin{proof}
Write \(2v=(m_1,m_2,m_3)\), where the integers have the same parity. Every path of \(n\) steps satisfies \(n\geq\max_i|m_i|\). Conversely, take \(n=\max_i|m_i|\). Each \(m_i\) is a sum of \(n\) signs; choosing these signs independently in the three coordinates gives a path of \(n\) allowed steps. Summation along this path proves \eqref{eq:displacement-distance}, from which the remaining assertions follow.
\end{proof}

\begin{proof}[Proof of Lemma~\ref{lem:local-cell}]
Write \(\eta=\eta(u)\). The half-diagonal differences have length at most \(\eta\), and the axis differences have length at most \(2\eta\). Thus
\begin{equation}\label{eq:local-array}
 \|d_p(u)\|^2\leq32\eta^2.
\end{equation}
We first show that the fourteen local inequalities exclude every remote generator. Define the reference gauge
\[
 t(y)=\max_{r\in\mathcal N}\frac{2r\cdot y}{|r|^2}.
\]
For \(y\ne0\), \(y/t(y)\in V\), so \(|y|\leq R_0t(y)\). If \(y\in K(d_p(u))\), the perturbed inequality for a maximizing \(r\) gives
\[
 t(y)\leq(1+2\eta)^2+4R_0\eta\,t(y).
\]
The half-diagonal inequalities give smaller coefficients than the axis inequalities. Hence
\begin{equation}\label{eq:local-radius}
 K(d_p(u))\subset B_{R_\eta}(0),
 \qquad
 R_\eta=R_0\frac{(1+2\eta)^2}{1-\sqrt5\eta}.
\end{equation}
For \(\eta\leq1/100\),
\[
 R_\eta<\frac35,\qquad (1-2\eta)\sqrt2>\frac65.
\]
Every lattice vector outside \(\mathcal N\cup\{0\}\) has length at least \(\sqrt2\). Its deformed vector \(w\) therefore satisfies \(|w|>2R_\eta\), by \eqref{eq:bilipschitz}. On the whole intersection,
\[
 w\cdot y\leq |w|R_\eta<|w|^2/2.
\]
Thus every omitted inequality is strict on the local intersection. This identifies the intersection with the actual Voronoi cell and also excludes remote contacts at vertices.

For analyticity, allow the fourteen vectors \(d_r\) to vary independently. The reference polytope is simple: at each vertex the three active normals have determinant of absolute value \(1/2\), and every inactive constraint has slack at least \(1/4\). Signed permutations reduce this check to \((1/2,1/4,0)\). Thus each vertex continues analytically, with its inactive inequalities strict.

The axis constraints keep all feasible points in a fixed bounded set for small \(d\): if \(a=\max_r|d_r|\) and \(b=\max_i|y_i|\), then
\[
 b\leq(1+a)^2/2+\sqrt3\,ab.
\]
Every other triple of labels has no feasible common point at the reference configuration. On the fixed bounded set, feasibility of each such triple is excluded by a positive separation; finiteness of the triples makes this exclusion uniform under small perturbations. Hence the continued vertices give the complete, unchanged face incidences.

Triangulate each face using the average of its vertices, and cone these triangles to the origin. This gives a fixed finite collection of nondegenerate tetrahedra with analytic vertices. Their quadratic moments, and therefore \(\Phi\), are analytic. Bounds for its third and fourth derivatives follow on a smaller closed parameter ball.

Finally, take the full cells corresponding to representatives of \(\Lambda/\Gamma\) and remove their boundary overlaps. They form a measurable fundamental region of \(\Gamma\). By periodicity, the energy integral over this region equals that over \(\Omega_\Gamma\). Translating each cell to its generator proves \eqref{eq:local-energy}. All labeled incidences are retained on short quotients.
\end{proof}

\begin{proof}[Proof of Proposition~\ref{prop:bond}]
We compute the variation of each interface, keeping track of the normalization; see also \cite{Urschel,Ruscher}. In every face sum below, we take one representative of each \(\Gamma\)-orbit of geometric unoriented faces in the infinite tessellation. For a face between \(p,q\), set
\[
 b_{pq}(x)=(x-p)\cdot u_p-(x-q)\cdot u_q,
 \qquad d=|q-p|.
\]
The normal velocity of its bisector, directed from the \(p\)-cell to the \(q\)-cell, is \(b_{pq}(x)/d\). The first variation is
\[
 DE_\Gamma(0)[u]=2\sum_p\int_{V_p}(p-x)\cdot u_p\,dx=0,
\]
because the reference cells are centered at their generators. The first boundary terms cancel since the two squared distances agree on their interface. Differentiating again gives
\begin{equation}\label{eq:interface-hessian}
 D^2E_\Gamma(0)[u,u]
 =2|V|\sum_p|u_p|^2
 -2\sum_{\{p,q\}}\frac1d\int_{F_{pq}}b_{pq}(x)^2\,dS.
\end{equation}
Writing \(x=(p+q)/2+y\), we have
\[
 b_{pq}(x)=\frac d2\,n\cdot(u_p+u_q)+y\cdot(u_p-u_q).
\]
The cross term integrates to zero by central symmetry of the face. The divergence theorem gives
\[
 |V|I=\sum_{F_r\subset\partial V}\frac{Ad}{2}nn^T.
\]
Use \((a+b)^2=2a^2+2b^2-(a-b)^2\) in \eqref{eq:interface-hessian}. The resulting coefficient for an unoriented bond is
\[
 B_r=\frac{Ad}{2}nn^T-\frac2d\int_{F_r}yy^T\,dS.
\]
Substitution of the moments in Lemma~\ref{lem:face-moments} gives \eqref{eq:bond-matrices}. Since no nonzero translation fixes a bounded face, the divergence identity retains both incidences of every face orbit. The formula therefore also applies to labeled loops and repeated endpoints on small quotients.
\end{proof}

Although the individual matrices \(B_r\) are indefinite, the bond sum is coercive. In the next section, we prove this by an exact estimate for the full Fourier symbol.

\section{Sharp coercivity and nonlinear stability}\label{sec:fourier}\label{sec:nonlinear}

We now prove Theorems~\ref{thm:stability} and~\ref{thm:strain}. We first establish sharp Fourier coercivity, classify its equality modes, and prove cubic cancellation on their span. The local cell representation then gives uniform nonlinear stability and allows simultaneous volume-preserving strain. The Fourier argument begins with the following matrix inequality.

\begin{lemma}[Sharp trigonometric matrix inequality]\label{lem:matrix}
For $x\in\R^3$, set $s_i=\sin x_i$, $c_i=\cos x_i$, $C=c_1c_2c_3$ and
$S=s_1^2+s_2^2+s_3^2$. The symmetric matrix
\begin{equation}\label{eq:symbol}
 M_{ii}(x)=26(1-C)-S+25s_i^2,
 \qquad M_{ij}(x)=41s_is_jc_\ell,
 \quad \{i,j,\ell\}=\{1,2,3\},
\end{equation}
satisfies $M(x)\succeq8(1-C)I$. The coefficient $8$ is optimal.
\end{lemma}

The proof reduces the determinant to two polynomial identities, whose exact coefficient formulas are recorded in Appendix~\ref{app:certificates}. The same identities will identify the equality modes.

\begin{proof}
Put
\[
 a=|c_1|,\qquad b=|c_2|,\qquad c=|c_3|,
 \qquad u=1-a^2,\quad v=1-b^2,\quad w=1-c^2.
\]
Diagonal sign conjugation removes the sine signs and makes the three off-diagonal cosine factors have the common sign $\sigma=\operatorname{sgn}C$: their product is $\sigma$, so the required conjugation exists. Zero sine or cosine coordinates follow by continuity. Thus $M-8(1-C)I$ is orthogonally equivalent to the symmetric matrix $N_\sigma$ with diagonal entries
\[
 A+25u,\quad A+25v,\quad A+25w,
 \qquad A=18(1-\sigma abc)-(u+v+w),
\]
and whose off-diagonal entries are
$\sigma41\sqrt{uv}\,c$, $\sigma41\sqrt{uw}\,b$ and
$\sigma41\sqrt{vw}\,a$. Writing $S=u+v+w$ and $\nu=\sigma abc$, direct expansion gives
\begin{equation}\label{eq:determinant}
\begin{split}
 F_\sigma(a,b,c):=\det N_\sigma
 ={}&A^3+25A^2S-1056A(uv+uw+vw)\\
 &+uvw(-110450+5043A+42025S+137842\nu).
\end{split}
\end{equation}
For the negative branch, the exact identity
\begin{equation}\label{eq:negative-certificate}
 F_-(a,b,c)=\sum_{i,j,k=0}^6
 Z^-_{ijk}\,a^i(1-a)^{6-i}b^j(1-b)^{6-j}c^k(1-c)^{6-k}
\end{equation}
has nonnegative integer coefficients and $Z^-_{000}=64000$.
Lemma~\ref{lem:negative-certificate-data} verifies the coefficients, with no implicit binomial factors in this basis. In particular,
$F_-\ge64000(1-a)^6(1-b)^6(1-c)^6>0$ on $[0,1)^3$.

For the positive branch, symmetry reduces the problem to
$0\le a\le b\le c\le1$. Introduce nonnegative variables $t,h,k,l$ and
$L=t+h+k+l$. Homogenizing to degree nine gives
\begin{equation}\label{eq:positive-certificate}
\begin{split}
 L^9F_+\left(\frac{t}{L},\frac{t+h}{L},\frac{t+h+k}{L}\right)
 &=6h^2(t+h+k)^3A_2B_2+lB(t,h,k,l),\\
 A_2&=17h^2+32hk+16ht+16k^2+32kt,\\
 B_2&=119h^2+246hk+299ht+123k^2+246kt+176t^2.
\end{split}
\end{equation}
The homogeneous polynomial $B$ has degree eight and exactly $161$ nonzero monomial coefficients, all positive; its $l^8$ coefficient is $64000$.
Lemma~\ref{lem:positive-certificate-data} gives an exact coefficient formula and proves these assertions.
Taking $(t,h,k,l)=(a,b-a,c-b,1-c)$ gives
$F_+(a,b,c)\ge64000(1-c)^9>0$ for $c<1$.
Each branch therefore has positive determinant on the connected region
$[0,1)^3$, using permutation symmetry for the positive branch. Since
$N_\sigma(0,0,0)=40I$, no eigenvalue can change sign there.
Both matrices are positive definite in this region and positive semidefinite on the closed cube by continuity.

Finally, for $x=(\theta,\theta,0)$ and $r=\sin^2\theta$,
\[
 M(x)=r\begin{pmatrix}49&41&0\\41&49&0\\0&0&24\end{pmatrix},
 \qquad M(x)(1,-1,0)^T=8r(1,-1,0)^T,
 \qquad 1-C=r.
\]
Thus no coefficient larger than $8$ is possible.
\end{proof}

We use the unitary Fourier transform on the finite abelian group $\Lambda/\Gamma$.
Its characters are $p\mapsto e^{ik\cdot p}$ with
$k\in\Gamma^*/\Lambda^*$, where reciprocal lattices include the factor $2\pi$.
By the bond formula \eqref{eq:bond}, direct summation gives the Hessian symbol
\begin{equation}\label{eq:bond-symbol}
 H(k)=2\sum_{r\in\mathcal R\cup\{e_1,e_2,e_3\}}
 (1-\cos(k\cdot r))B_r=\frac1{96}M(k/2),
 \qquad
 \sum_{r\in\mathcal R}2(1-\cos(k\cdot r))=8(1-C).
\end{equation}
Lemma~\ref{lem:matrix} and Parseval's identity consequently yield
\begin{equation}\label{eq:sharp-hessian}
 D^2E_\Gamma(0)[v,v]\ge\frac1{96}D_\Gamma(v).
\end{equation}
The constant is attained on a finite quotient. Indeed, take $\Gamma=2\Z^3$ and
\[
 v_p=(-1)^{p_1+p_2}(1,-1,0).
\]
Here $p_1+p_2\in\Z$ on $\Lambda$, the wave vector is $(\pi,\pi,0)$,
and the ratio in \eqref{eq:sharp-hessian} is exactly $1/96$.
Since the half-diagonal steps generate $\Lambda$, the nullspace of $D_\Gamma$ consists exactly of translations. We now determine the full equality set in \eqref{eq:sharp-hessian}.

\begin{proposition}[Complete quadratic equality set]\label{prop:equality}
Let $N(x)=M(x)-8(1-C)I$. Apart from the translation frequency, $N(x)$ is singular exactly when, after a coordinate permutation,
\[
 \cos x_3=\varepsilon\in\{-1,1\},
 \qquad \cos x_2=\varepsilon\cos x_1.
\]
With $r=\sin^2x_1>0$, its kernel is
\[
 \ker N(x)=\Span\{(\sin x_1,-\varepsilon\sin x_2,0)\}.
\]
Equality in \eqref{eq:sharp-hessian} holds if and only if every nontranslation Fourier coefficient belongs to this kernel at its frequency; an arbitrary translation may be added.
\end{proposition}

\begin{proof}
Strict positivity in the open cube leaves only boundary points, where we may take $c=1$. There \eqref{eq:determinant} factors as
\[
\begin{split}
 F_+(a,b,1)&=6(a-b)^2(16+a^2+b^2-18ab)
                 (123-4a^2+61ab-4b^2),\\
 F_-(a,b,1)&=6(a+b)^2(16+a^2+b^2+18ab)
                 (123-4a^2-61ab-4b^2).
\end{split}
\]
The final factors are at least $115$ and $54$, while
$16+a^2+b^2-18ab=16(1-ab)+(a-b)^2$ and
$16+a^2+b^2+18ab\ge16$.
Thus the positive branch is singular exactly when $a=b$, and the negative branch only when $a=b=0$. Restoring signs gives precisely the stated cosine conditions, under which
\[
 N(x)=\begin{pmatrix}
 41r&41\varepsilon\sin x_1\sin x_2&0\\
 41\varepsilon\sin x_1\sin x_2&41r&0\\
 0&0&16r
 \end{pmatrix},
\]
which proves the kernel formula.
For $r=0$, these conditions give $C=1$ and $M=0$, precisely the translation frequency modulo $\Lambda^*$. If all sines vanish and $C=-1$, instead $N=36I$.
The modal equality characterization now follows from Parseval's identity.
\end{proof}

The equality fields are null directions of the defect
$D^2E_\Gamma(0)-D_\Gamma/96$. Every nontranslation equality field still has positive second variation.
The scaled reciprocal lattice is
\[
 \frac1{2\pi}\Lambda^*=D_3
 :=\{m\in\Z^3:m_1+m_2+m_3\ \text{is even}\}.
\]
Modulo these reciprocal shifts, every nontranslation equality frequency
lies on one of the six unoriented lines $k=q\ell$, with $q$ taken modulo
$2\pi$. The roots and their polarizations are given uniformly by
\begin{equation}\label{eq:root-polarizations}
 \ell=e_i+\varepsilon e_j,\qquad
 v_\ell=e_i-\varepsilon e_j,
 \qquad 1\leq i<j\leq3,\quad \varepsilon\in\{-1,1\}.
\end{equation}

\begin{proposition}[Cubic cancellation on the entire equality space]\label{prop:cubic}
If a real periodic displacement $v$ satisfies equality in
\eqref{eq:sharp-hessian}, then
\[
 D^3E_\Gamma(0)[v,v,v]=0.
\]
Equivalently, the complex trilinear third variation vanishes whenever all three arguments belong to the reference equality space.
\end{proposition}

\begin{proof}
By Lemma~\ref{lem:local-cell}, the third variation is well defined.
Let $G$ be the signed coordinate permutation group. Passing to the common period
$\Gamma'=\bigcap_{R\in G}R\Gamma$ multiplies the energy, $D_\Gamma$, and all variations by $[\Gamma:\Gamma']$, preserving equality. We may therefore assume that the quotient is $G$-invariant and write $T$ for the complex trilinear third variation.

Translation invariance gives zero whenever an argument is constant, or when
$k_1+k_2+k_3\notin2\pi D_3$. It remains to consider resonant triples of nontranslation equality modes.

An independent root triple is an integer basis of $D_3$: its determinant is a nonzero even integer bounded by $(\sqrt2)^3<3$, hence has absolute value two. Resonance then forces all $q_i\in2\pi\Z$, which gives only constant arguments.

Similarly, every nonparallel root pair generates its plane's intersection with $D_3$. Up to signs and coordinate permutation, the pair is
$(1,1,0),(1,-1,0)$ or $(1,1,0),(1,0,1)$.
The intersections consist respectively of
$(s+t,s-t,0)$ and $(b+c,b,c)$, with all displayed parameters integral.
Thus, for a triple with exactly two root directions, resonance makes the unpaired argument constant.

If all three directions agree, a signed coordinate interchange fixes the frequency line and reverses its polarization. For $(1,1,0)$ this is the interchange of the first two coordinates. Point-group invariance and trilinearity give $T=-T=0$.

For three distinct coplanar directions, an orthogonal pair has no third root in its plane; every remaining triple reduces, up to signs and signed coordinate permutation, to
\[
 \ell_1=(1,1,0),\qquad \ell_2=(1,0,1),\qquad
 \ell_3=(0,1,-1),\qquad \ell_1=\ell_2+\ell_3.
\]
The pair calculation gives
$(q_1,q_2,q_3)\equiv(q,-q,-q)\pmod{2\pi}$.
Choose polarizations $v_1=(1,-1,0)$, $v_2=(1,0,-1)$ and $v_3=(0,1,1)$.
Interchanging the second and third coordinates sends
\[
 ((k_1,v_1),(k_2,v_2),(k_3,v_3))
 \longmapsto
 ((-k_2,v_2),(-k_1,v_1),(-k_3,v_3)).
\]
Hence point-group invariance and symmetry identify the coupling with that at the three reversed frequencies. Spatial inversion reverses both frequencies and polarizations, so trilinearity gives the opposite sign. The coupling vanishes.
All root triples are covered, including boundary frequencies modulo $2\pi D_3$. Fourier expansion proves the cubic identity, and polarization proves its trilinear form.
\end{proof}

The sharp reference estimate and its equality space are now established. We next use the local cell representation from Section~\ref{sec:geometry} to compare the Hessian at a perturbed configuration with the reference Hessian. This gives uniform nonlinear stability and the Taylor bounds needed for the refined loss estimates.

\begin{proposition}[Uniform Hessian comparison]\label{prop:hessian-comparison}
For \(\eta(u)\leq\min\{1/100,\delta/(4\sqrt2)\}\),
\begin{equation}\label{eq:hessian-comparison}
 \left|D^2E_\Gamma(u)[v,v]-D^2E_\Gamma(0)[v,v]\right|
 \leq20\sqrt2M_3\eta(u)D_\Gamma(v).
\end{equation}
The third-order remainder satisfies
\begin{equation}\label{eq:cubic-remainder}
 \left|E_\Gamma(u)-E_\Gamma(0)-\tfrac12D^2E_\Gamma(0)[u,u]\right|
 \leq\frac{10\sqrt2M_3}{3}\eta(u)D_\Gamma(u).
\end{equation}
On the same neighborhood,
\begin{equation}\label{eq:taylor}
 \begin{split}
 E_\Gamma(u)-E_\Gamma(0)
 &=\tfrac12D^2E_\Gamma(0)[u,u]
 +\tfrac16D^3E_\Gamma(0)[u,u,u]+R_4(u),\\
 |R_4(u)|&\leq K_4\eta(u)^2D_\Gamma(u),
 \qquad K_4=\frac{20M_4}{3}.
 \end{split}
\end{equation}
\end{proposition}

By the local cell representation, it suffices to control the fourteen-neighbor difference sum in terms of \(D_\Gamma\). The cell derivative bounds then give both the Hessian comparison and the Taylor remainders.

\begin{proof}
The positive-axis sum has symbol \(4\sum_i\sin^2x_i\), whereas \(D_\Gamma\) has symbol \(8(1-C)\). The arithmetic--geometric mean inequality gives, for every sign of \(C\),
\[
 \sum_i\sin^2x_i\leq3(1-|C|^{2/3})\leq3(1-C).
\]
Parseval's identity therefore yields
\[
 \sum_{p\in\Lambda/\Gamma}\sum_{i=1}^3|v_{p+e_i}-v_p|^2
 \leq\frac32D_\Gamma(v).
\]
Including both signs gives
\begin{equation}\label{eq:neighbor-comparison}
 \sum_p\|d_p(v)\|^2
 =\sum_p\sum_{r\in\mathcal N}|v_{p+r}-v_p|^2
 \leq5D_\Gamma(v).
\end{equation}
The constant is sharp, for example when all three cosine coordinates vanish.

Apply the mean-value estimate for \(D^2\Phi\) in \eqref{eq:local-energy}. Using \eqref{eq:local-array} and \eqref{eq:neighbor-comparison}, we obtain
\[
 \begin{split}
 \left|D^2E_\Gamma(u)[v,v]-D^2E_\Gamma(0)[v,v]\right|
 &\leq M_3\max_p\|d_p(u)\|\sum_p\|d_p(v)\|^2\\
 &\leq20\sqrt2M_3\eta(u)D_\Gamma(v).
 \end{split}
\]
Integration along \(tu\) with weight \(1-t\) gives \eqref{eq:cubic-remainder}, since \(\int_0^1t(1-t)\,dt=1/6\). The fourth-order cell remainder is at most
\[
 \frac{M_4}{24}\max_p\|d_p(u)\|^2\sum_p\|d_p(u)\|^2
 \leq\frac{20M_4}{3}\eta(u)^2D_\Gamma(u),
\]
which proves \eqref{eq:taylor}.
\end{proof}

\begin{proof}[Proof of Theorem~\ref{thm:stability}]
Set
\begin{equation}\label{eq:eta0}
 \eta_0=\min\left\{\frac1{100},\frac{\delta}{4\sqrt2},
 \frac1{3840\sqrt2M_3}\right\}.
\end{equation}
The sharp reference estimate and \eqref{eq:hessian-comparison} give
\[
 D^2E_\Gamma(u)[v,v]\geq\frac1{192}D_\Gamma(v),
 \qquad \eta(u)\leq\eta_0.
\]
Since \(DE_\Gamma(0)=0\), integration along \(tu\) yields
\[
 E_\Gamma(u)-E_\Gamma(0)
 =\int_0^1(1-t)D^2E_\Gamma(tu)[u,u]\,dt
 \geq\frac1{384}D_\Gamma(u).
\]
The nearest-neighbor directions generate \(\Lambda\), so equality in energy forces \(u\) to be constant. If \(u\) is stationary, the same argument gives
\[
 0=DE_\Gamma(u)[u]-DE_\Gamma(0)[u]
 =\int_0^1D^2E_\Gamma(tu)[u,u]\,dt
 \geq\frac1{192}D_\Gamma(u).
\]
Thus every stationary configuration in this neighborhood is a translate of BCC.
\end{proof}

All constants in \eqref{eq:eta0} depend only on the reference cell.

The hypothesis on displacement differences also permits large long-wave displacements. For instance, when \(4\mid m\), the field
\[
 u_p=\sqrt m\,\sin(2\pi p_1/m)e_2,\qquad\Gamma=m\Z^3,
\]
has \(\eta(u)\leq\pi/\sqrt m\to0\), while its amplitude modulo translations is \(\sqrt m\to\infty\). A period-independent coercivity estimate in the displacement \(\ell^2\) norm is impossible, since the long-wave eigenvalues tend to zero. A fixed-period Poincar\'e inequality supplies that norm with a period-dependent constant.

\begin{corollary}[Compactly supported perturbations]\label{cor:compact}
For a compactly supported displacement of the infinite lattice with \(\eta(u)\leq\eta_0\), the excess
\[
 \Delta E(u)=\int_{\R^3}
 \bigl(\dist(x,\Lambda_u)^2-\dist(x,\Lambda)^2\bigr)\,dx
\]
is finite and satisfies \(\Delta E(u)\geq D(u)/384\), where \(D(u)\) is the full infinite-lattice bond sum.
\end{corollary}

\begin{proof}
By locality, a cell remains unchanged unless its generator or one of its fourteen neighbors belongs to the finite support. The changed cells are bounded, so the integrand vanishes outside a bounded set. We embed these cells and the displacement support in a sufficiently large period box, with a buffer exceeding the finite interaction range. The periodic extension has the same excess and bond sum per box, and the result follows from Theorem~\ref{thm:stability}.
\end{proof}

We next turn to homogeneous strain. The proof compares periodic displacements at a fixed strained lattice and then uses the following exact second variation of the homogeneous lattice energy.

Agrell and Allen \cite{AgrellAllen} proved that local optimality among lattice quantizers forces the second-moment tensor to be isotropic. Formula \eqref{eq:cell-moment} verifies this necessary condition for BCC. To obtain the quantitative strain estimate, we compute the full elastic quadratic form.

\begin{proposition}[Exact elastic quadratic form]\label{prop:elastic}
Let \(S=S^T\), \(\tr S=0\), and write \(W(F)=\E_\Gamma(F,0)/N\). Then \(W\) is independent of the chosen period, and
\begin{equation}\label{eq:elastic-hessian}
 \left.\frac{d^2}{dt^2}W(e^{tS})\right|_{t=0}
 =\frac1{48}\sum_iS_{ii}^2+\frac18\sum_{i<j}S_{ij}^2.
\end{equation}
\end{proposition}

\begin{proof}
Pull back by \(F=e^{tS}\), whose determinant is one. On the fixed reference torus the integrand is the minimum of
\[
 (x-p)^TQ(t)(x-p),\qquad Q(t)=e^{2tS}.
\]
By \eqref{eq:cell-moment}, the first variation vanishes and the bulk second derivative is
\[
 4\int_Vx^TS^2x\,dx=\frac{19}{192}\tr S^2.
\]
On a face \(r=dn\), write \(x=r/2+y\). The first metric variations of the two branches differ by \(4y\cdot Sr\), so the interface correction per unoriented face is
\[
 -\frac8d\int_{F_r}(y\cdot Sr)^2\,dS
 =-8A\beta d\bigl(|Sn|^2-(n^TSn)^2\bigr).
\]
The face constants are \(8A\beta d=1/96\) on an axis face and \(15/256\) on a half-diagonal face. Since \(\tr S=0\),
\[
 \begin{aligned}
 \sum_i\bigl(|Se_i|^2-S_{ii}^2\bigr)&=2\sum_{i<j}S_{ij}^2,\\
 \sum_{r\in\mathcal R}\bigl(|Sn_r|^2-(n_r^TSn_r)^2\bigr)
 &=\frac43\sum_iS_{ii}^2+\frac89\sum_{i<j}S_{ij}^2.
 \end{aligned}
\]
Subtracting the corresponding face corrections from the bulk term gives \eqref{eq:elastic-hessian}.
\end{proof}

\begin{proof}[Proof of Theorem~\ref{thm:strain}]
Rotational invariance reduces the proof to \(F=e^S\). The deformed relative vectors have length at least
\[
 \bigl(1-\|F-I\|-2\eta(u)\bigr)|q-p|.
\]
For \(\|F-I\|+\eta(u)\) small, the local radius estimate and the strict separation from remote generators in Lemma~\ref{lem:local-cell} therefore persist. Hence
\[
 \E_\Gamma(F,u)=\sum_{p\in\Lambda/\Gamma}
 \Phi\bigl(((F-I)r+u_{p+r}-u_p)_{r\in\mathcal N}\bigr).
\]
The same mean-value estimate as in Proposition~\ref{prop:hessian-comparison} gives
\[
 \left|D_u^2\E_\Gamma(F,u)[v,v]-D^2E_\Gamma(0)[v,v]\right|
 \leq C\bigl(\|F-I\|+\eta(u)\bigr)D_\Gamma(v),
\]
with \(C\) independent of the period. Choose the radius so that
\[
 D_u^2\E_\Gamma(F,u)[v,v]\geq D_\Gamma(v)/192.
\]
Central symmetry of the cells of \(F\Lambda\) gives \(D_u\E_\Gamma(F,0)=0\). Integration along \(tu\) yields
\begin{equation}\label{eq:strained-split}
 \E_\Gamma(F,u)-\E_\Gamma(F,0)\geq D_\Gamma(u)/384.
\end{equation}
Analyticity and Proposition~\ref{prop:elastic} give \eqref{eq:strain-expansion}. Its quadratic term is at least \(|S|_F^2/96\); after absorbing the cubic remainder,
\[
 W(e^S)-W(I)\geq |S|_F^2/192.
\]
Together with \eqref{eq:strained-split}, this proves the estimate and its equality assertion.
\end{proof}

\section{Nonlinear loss and sharp period dependence}\label{sec:fixed}\label{sec:twoscale}
We now prove Theorems~\ref{thm:loss} and~\ref{thm:crossover}. For a fixed period, cubic cancellation and the transverse spectral gap give a quadratic correction to the coefficient $1/192$. An exact four-phase construction shows that this loss order is attained on $4\Z^3$. We then track the spectral gap as the period grows and use a nearby-frequency coupling to obtain the sharp crossover and the all-period linear loss.

On the mean-zero displacement space, use the inner product associated with $D_\Gamma$, and put
\[
 g_\Gamma(u)=D^2E_\Gamma(0)[u,u]-\frac1{96}D_\Gamma(u),
 \qquad K_\Gamma=\ker g_\Gamma.
\]
If $K_\Gamma^\perp\ne\{0\}$, define
\[
 \beta_\Gamma=
 \min_{0\ne w\in K_\Gamma^\perp}
 \frac{g_\Gamma(w)}{D_\Gamma(w)}>0.
\]
The sharp reference inequality and finite dimensionality imply that $\beta_\Gamma>0$. We use the fourth-derivative bound \(M_4\) from Lemma~\ref{lem:local-cell}.

\begin{proposition}[Quadratic loss for a fixed period]\label{prop:fixed-loss}
For every finite-index period sublattice $\Gamma\subset\Lambda$, there is $C_\Gamma<\infty$ such that
\begin{equation}\label{eq:fixed-loss}
 E_\Gamma(u)-E_\Gamma(0)
 \geq \left(\frac1{192}-C_\Gamma\eta(u)^2\right)D_\Gamma(u)
\end{equation}
whenever $\eta(u)$ is sufficiently small. One may take
\[
 C_\Gamma=\frac{125NM_3^2}{2\beta_\Gamma}+\frac{20M_4}{3},
\]
with the first term omitted when $K_\Gamma^\perp=\{0\}$.
\end{proposition}
\begin{proof}
Subtract the mean of $u$ and take the $D_\Gamma$-orthogonal decomposition
$u=v+w$, with $v\in K_\Gamma$ and $w\in K_\Gamma^\perp$. Write
$\|a\|_D=D_\Gamma(a)^{1/2}$ and $T=D^3E_\Gamma(0)$. Then
\[
 g_\Gamma(u)=g_\Gamma(w)\geq\beta_\Gamma\|w\|_D^2,
 \qquad \|v\|_D\leq\|u\|_D.
\]
The local derivative bound and~\eqref{eq:neighbor-comparison} give
\begin{equation}\label{eq:cubic-fixed-bound}
 |T[a,b,c]|\leq B\|a\|_D\|b\|_D\|c\|_D,
 \qquad B=5\sqrt5\,M_3.
\end{equation}
Indeed, the local cell formula bounds the left side by
$M_3\sum_p\|d_p(a)\|\|d_p(b)\|\|d_p(c)\|$, and the product of the
three $\ell^2$ norms bounds this sum. Each norm is at most $\sqrt5$
times the corresponding $D$ norm.
Since $T[v,v,v]=0$ by Proposition~\ref{prop:cubic}, telescoping the
three tensor slots and using $D_\Gamma(u)\leq4N\eta(u)^2$ yields
\[
 \frac16|T[u,u,u]|
 \leq\frac B2\|u\|_D^2\|w\|_D
 \leq B\sqrt{\frac N{\beta_\Gamma}}\,
 \eta(u)\|u\|_D\sqrt{g_\Gamma(u)}.
\]
Moreover, the fourth-order remainder satisfies
\[
 |R_4(u)|\leq\frac{M_4}{24}
 \max_p\|d_p(u)\|^2\sum_p\|d_p(u)\|^2
 \leq\frac{20M_4}{3}\eta(u)^2D_\Gamma(u),
\]
since the maximum is at most $32\eta(u)^2$ and the sum is at most
$5D_\Gamma(u)$. Taylor's formula and completion of the square now give
\begin{align*}
 E_\Gamma(u)-E_\Gamma(0)
 &\geq\frac1{192}D_\Gamma(u)+\frac12g_\Gamma(u)
 -B\sqrt{\frac N{\beta_\Gamma}}\eta(u)\|u\|_D\sqrt{g_\Gamma(u)}
 -\frac{20M_4}{3}\eta(u)^2D_\Gamma(u)\\
 &\geq\left(\frac1{192}
 -\frac{B^2N}{2\beta_\Gamma}\eta(u)^2
 -\frac{20M_4}{3}\eta(u)^2\right)D_\Gamma(u).
\end{align*}
If $K_\Gamma^\perp=\{0\}$, both $g_\Gamma$ and the cubic term vanish,
and the remainder estimate proves the assertion without the first term
in $C_\Gamma$.
\end{proof}

\begin{proposition}[Nonlinear obstruction to the sharp quadratic coefficient]\label{prop:four-phase}
For every $\rho>0$, there is a $4\Z^3$-periodic displacement $u$ with
$0<\eta(u)<\rho$ such that
\[
 E_{4\Z^3}(u)-E_{4\Z^3}(0)<\frac1{192}D_{4\Z^3}(u).
\]
\end{proposition}
\begin{proof}
Put $\Gamma=4\Z^3$, so $N=128$. For $a=p_1+p_2\in\Z$, define
\begin{equation}\label{eq:four-phase-field}
 \begin{aligned}
 w_p&=\bigl(\sin(\pi a/2)+\cos(\pi a/2)\bigr)(1,-1,0),\\
 z_p&=(-1)^a(1,1,0),\qquad
 u_p(t)=tw_p+\lambda t^2z_p,\qquad\lambda=\frac4{41}.
 \end{aligned}
\end{equation}
The scalar phase factors are $(1,1,-1,-1)$ for $w$ and $(1,-1,1,-1)$ for $z$, with $32$ sites in each phase. The field $w$ is an equality mode, and counting the four nearest-neighbor bonds gives
\begin{equation}\label{eq:fourphase-D}
 \frac{D_\Gamma(u(t))}{N}=8t^2+16\lambda^2t^4,
 \qquad \eta(u(t))^2=8t^2+8\lambda^2t^4.
\end{equation}
The exact local cell calculation in Appendix~\ref{app:cells} gives
\begin{equation}\label{eq:fourphase-energy}
 \frac{E_\Gamma(u(t))}{N}=J_\lambda(t)
 =\frac{57+32t^2-(192\lambda^2+128\lambda)t^4-12288\lambda^4t^8}
 {768(1-16\lambda^2t^4)}.
\end{equation}
For every fixed $\lambda$, the identity holds for the actual Voronoi cells in a neighborhood of $t=0$. Substituting $\lambda=4/41$ gives
\begin{equation}\label{eq:fourphase-defect}
 \frac{E_\Gamma(u(t))-E_\Gamma(0)-D_\Gamma(u(t))/192}{N}
 =-\frac{t^4(68921-53792t^2+11264t^4)}{5043(1681-256t^4)}.
\end{equation}
The right side is strictly negative for all sufficiently small nonzero $t$, and $\eta(u(t))\to0$. This proves the assertion.
\end{proof}

To explain the choice of the second-order correction, expand~\eqref{eq:fourphase-energy} for a general fixed $\lambda$:
\[
 \frac{E_\Gamma(u(t))-E_\Gamma(0)-D_\Gamma(u(t))/192}{N}
 =\left(-\frac\lambda6+\frac{41\lambda^2}{48}\right)t^4+O(t^6).
\]
The coefficient is minimized at $\lambda=4/41$. Hence the correction in the polarization $(1,1,0)$ lowers the energy ratio of a leading displacement in the equality space. This effect is not detected by straight paths in that space. By~\eqref{eq:fourphase-D}, with $\rho^2=\eta(u(t))^2$, we obtain
\begin{equation}\label{eq:four-phase-ratio}
 \frac{E_\Gamma(u(t))-E_\Gamma(0)}{D_\Gamma(u(t))}
 =\frac1{192}-\frac{t^2}{984}+O(t^4)
 =\frac1{192}-\frac{\rho^2}{7872}+O(\rho^4).
\end{equation}
Every sufficiently small positive $\rho$ is realized by a unique small positive $t$.

\begin{corollary}[Optimal loss order on one fixed quotient]\label{cor:fixed-optimal}
For $\Gamma=4\Z^3$, set
\[
 c_4(\rho)=\inf_{0<\eta(u)\leq\rho}
 \frac{E_\Gamma(u)-E_\Gamma(0)}{D_\Gamma(u)}.
\]
For all sufficiently small $\rho>0$,
\[
 \frac{\rho^2}{15744}\leq\frac1{192}-c_4(\rho)
 \leq C_{4\Z^3}\rho^2.
\]
\end{corollary}
\begin{proof}
The upper bound follows from Proposition~\ref{prop:fixed-loss}, and the lower bound follows from~\eqref{eq:four-phase-ratio}. These bounds determine the order of loss; they do not determine its optimal leading constant.
\end{proof}

The fixed-period estimate depends on the transverse gap $\beta_\Gamma$. We now determine how this gap changes along growing cubic periods and how it governs the transition from quadratic to linear loss. The theorem below shows that the quadratic-loss constants are at least of order $m^2$, even when their neighborhoods shrink arbitrarily fast with $m$.

For $\Gamma_m=m\Z^3$, write $N_m=2m^3$ and abbreviate
\[
 E_m=E_{\Gamma_m},\quad D_m=D_{\Gamma_m},\quad
 g_m=g_{\Gamma_m},\quad K_m=K_{\Gamma_m},\quad
 \beta_m=\beta_{\Gamma_m},\quad T_m=D^3E_m(0).
\]
We work on mean-zero fields, since translations do not affect these quantities. The defect and bond symbols are respectively $N(x)/96$ and $8(1-c_1c_2c_3)I$, where
\[
 N(x)=M(x)-8(1-c_1c_2c_3)I,\qquad x=k/2,
 \qquad c_i=\cos x_i,\quad s_i=\sin x_i.
\]

\begin{theorem}[Necessary growth on arbitrarily shrinking neighborhoods]\label{thm:necessary-growth}
Suppose that, for each $4\mid m$, a constant $A_m$ makes
\[
 E_m(u)-E_m(0)\geq
 \left(\frac1{192}-A_m\eta(u)^2\right)D_m(u)
\]
valid on some positive small-gradient neighborhood, whose radius may depend on $m$. Then
\begin{equation}\label{eq:necessary-growth}
 \liminf_{\substack{m\to\infty\\4\mid m}}
 \left(\frac{2\pi}{m}\right)^2A_m\geq\frac2{171}.
\end{equation}
In particular,
$\liminf_{m\to\infty,\,4\mid m}C_m^{\rm loc}/m^2\geq1/(342\pi^2)$.
\end{theorem}

The upper bounds follow from a transverse spectral gap and a bounded projection onto the equality space. For the reverse bounds, we couple two equality modes to a nearby frequency whose quadratic defect is of order $m^{-2}$.

\begin{lemma}[Quadratic transverse gap]\label{lem:gap}
Let $D_3=\{a\in\Z^3:a_1+a_2+a_3\text{ is even}\}$ and
$\mathbb T=\R^3/(\pi D_3)$. Let $\mathcal Z\subset\mathbb T$ be the union of the six circles
\[
 x=\theta\ell,\qquad
 \ell\in\mathcal L:=\{(1,1,0),(1,-1,0),(1,0,1),(1,0,-1),(0,1,1),(0,1,-1)\}.
\]
There are constants $c_*,b>0$, independent of $m$, such that
\begin{equation}\label{eq:gap}
 \lambda_{\min}N(x)\geq c_*\dist_{\mathbb T}(x,\mathcal Z)^2,
 \qquad \beta_m\geq \frac{b}{m^2}.
\end{equation}
\end{lemma}

\begin{proof}
Proposition~\ref{prop:equality} gives the zero set $\mathcal Z$. The root-pair saturation proved in Proposition~\ref{prop:cubic} shows that distinct circles meet only at translation: $\theta\ell-\theta'\ell'\in\pi D_3$ forces $\theta,\theta'\in\pi\Z$.

Away from translation, symmetry reduces the transverse calculation to $x=(\theta+h,\theta-h,z)$, with $r=\sin^2\theta>0$. At $h=z=0$, the eigenvalues of $N$ are $0,82r,16r$, with weak polarization $(1,-1,0)$. Differentiating the explicit symbol and using the second-order formula for a simple eigenvalue gives
\begin{equation}\label{eq:transverse-curvature}
 \lambda_{\min}N(\theta+h,\theta-h,z)
 =\frac{12(176-53r)}{41}h^2+\frac{16+23r}{2}z^2
 +O_\theta((|h|+|z|)^3).
\end{equation}
Both coefficients are positive, which proves the desired bound on compact portions of the circles away from translation.

Near translation, write $N(x)=N_2(x)+O(|x|^4)$, where
\[
 N_2(x)=8|x|^2I-16\diag(x_1^2,x_2^2,x_3^2)+41xx^T,
\]
and observe the exact identity
\begin{equation}\label{eq:origin-sos}
 v^TN_2(x)v
 =8\sum_{i<j}(x_iv_j+x_jv_i)^2+33(x\cdot v)^2.
\end{equation}
Its kernel is nonzero exactly on the six root lines, and is one-dimensional there. Indeed, the equations on the right force $v=0$ when $x$ has one or three nonzero coordinates; with two nonzero coordinates, a kernel exists exactly when their magnitudes agree.

For $x=t\omega$, $|\omega|=1$, the matrix $t^{-2}N(t\omega)$ extends analytically to $t=0$ with value $N_2(\omega)$. At $\omega=(1,1,0)/\sqrt2$, its simple lowest eigenvalue has positive quadratic angular coefficients $2112/41$ and $8$ in tangent coordinates $(h,-h,z)$; these follow directly from the displayed formula for $N_2$. Symmetry gives the same conclusion at every unit root direction. For small $t$, the eigenvalue is still zero at these exact directions, its angular gradient vanishes by nonnegativity, and its angular Hessian remains uniformly positive. Taylor's formula therefore bounds it below by the squared angular distance to the root directions. Away from them, \eqref{eq:origin-sos} and compactness give a positive lower bound. Multiplying by $t^2$ proves the first inequality near translation. Compactness on the rest of $\mathbb T$ completes it.

Allowed frequencies have $x=\pi n/m$ modulo $\pi D_3$. Their squared distance to a translated $(1,\pm1,0)$ root line is
\[
 (x_3-\pi a_3)^2+
 \tfrac12\bigl((x_1-\pi a_1)\mp(x_2-\pi a_2)\bigr)^2,
 \qquad a\in D_3.
\]
Thus every frequency outside $\mathcal Z$ has squared distance at least $\pi^2/(2m^2)$. Since the bond symbol is at most $16$, the quotient $g_m/D_m$ there is at least $c_*\pi^2/(3072m^2)$. At a nontranslation equality frequency, the complementary eigenvalues $82r,16r$ and bond symbol $8r$ give quotient at least $1/48$. Parseval proves the second inequality with $b=\min\{c_*\pi^2/3072,1/48\}$.
\end{proof}

\begin{lemma}[Uniform projection estimate]\label{lem:projection}
Let $P_m$ be the $D_m$-orthogonal projection onto $K_m$. Then
\[
 \eta(P_mu)\leq6\eta(u),\qquad D_m(P_mu)\leq D_m(u).
\]
Moreover, with the universal third-derivative bound $M_3\geq1$,
\begin{equation}\label{eq:cubic-bound}
 |T_m[u,u,u]|
 \leq160\sqrt2M_3\eta(u)\sqrt{D_m(u)}
             \sqrt{D_m((I-P_m)u)}.
\end{equation}
\end{lemma}

\begin{proof}
For $\ell\in\mathcal L$, average translations over
$H_\ell=\{p\in\Lambda/\Gamma_m:\ell\cdot p\equiv0\pmod m\}$, and denote this operator by $A_\ell$. It projects onto the corresponding equality circle, since $p\mapsto\ell\cdot p\bmod m$ maps onto $\Z/m\Z$. If $\Pi_\ell$ projects onto its weak polarization, the equality classification gives
\[
 P_m=\sum_{\ell\in\mathcal L}\Pi_\ell A_\ell
\]
on mean-zero fields. Each summand contracts every maximum bond-difference norm, so $\eta(P_mu)\leq6\eta(u)$. The scalar bond symbol also gives $D_m(P_mu)\leq D_m(u)$.

The local-array and neighbor estimates \eqref{eq:local-array}, \eqref{eq:neighbor-comparison}, together with Cauchy--Schwarz, give the polarized bound
\begin{equation}\label{eq:polarized-cubic-bound}
 |T_m[a,b,c]|
 \leq20\sqrt2M_3\eta(b)\sqrt{D_m(a)}\sqrt{D_m(c)}.
\end{equation}
Put $v=P_mu$, $w=u-v$. By Proposition~\ref{prop:cubic},
\[
 T_m[u,u,u]=T_m[w,u,u]+T_m[v,w,u]+T_m[v,v,w].
\]
Use symmetry in \eqref{eq:polarized-cubic-bound} to place $u,u,v$ in its middle slot, respectively. The projection bounds give the factor $1+1+6=8$, proving \eqref{eq:cubic-bound}.
\end{proof}

\begin{proof}[Proof of the upper estimates in Theorem~\ref{thm:crossover}]
By \eqref{eq:taylor},
\begin{equation}\label{eq:fourth-remainder}
 \begin{split}
 E_m(u)-E_m(0)-\frac{D_m(u)}{192}
 &=\frac12g_m(u)+\frac16T_m[u,u,u]+R_4(u),\\
 |R_4(u)|&\leq K_4\eta(u)^2D_m(u),\qquad K_4=\frac{20M_4}{3}.
 \end{split}
\end{equation}
Its radius is independent of $m$. Since $g_m((I-P_m)u)=g_m(u)$, Lemma~\ref{lem:projection} bounds the cubic term by
\[
 \frac16|T_m[u,u,u]|
 \leq\frac{80\sqrt2M_3}{3\sqrt{\beta_m}}\eta(u)
 \sqrt{D_m(u)g_m(u)}.
\]
Completing the square and applying Lemma~\ref{lem:gap} gives
\begin{equation}\label{eq:matching}
 \begin{split}
 E_m(u)-E_m(0)
 &\geq\left[\frac1{192}-
 \left(\frac{6400M_3^2}{9\beta_m}+\frac{20M_4}{3}\right)\eta(u)^2\right]D_m(u)\\
 &\geq\left(\frac1{192}-Cm^2\eta(u)^2\right)D_m(u).
 \end{split}
\end{equation}
Projected fields enter only derivatives at zero, so their size does not reduce the uniform radius. Combining \eqref{eq:matching} with \eqref{eq:cubic-remainder} proves
$1/192-c_m(\rho)\leq C\min\{m^2\rho^2,\rho\}$ and $C_m^{\rm loc}\leq Cm^2$.
\end{proof}

\begin{proposition}[Two-scale family]\label{prop:two-scale}
Let $4\mid m$, $m\geq8$, and put
\[
 \ell=(1,1,0),\quad v=(1,-1,0),\quad
 h=(1,0,1),\quad w=(1,0,-1),\quad
 \varepsilon=\frac{2\pi}{m},\quad q=\frac\pi2.
\]
Define
\begin{equation}\label{eq:two-scale}
 V_m(p)=v\{\cos\psi_p+\cos(\psi_p+\theta_p)\}
          +\frac{w}{\varepsilon}\sin\theta_p,
 \qquad \psi_p=q\ell\cdot p,\quad\theta_p=\varepsilon h\cdot p.
\end{equation}
This real $\Gamma_m$-periodic field satisfies
\begin{align}
 \eta(V_m)&\leq5\sqrt2,\label{eq:two-scale-eta}\\
 \frac{g_m(V_m)}{N_m}&=\frac{19}{128}(1-\cos\varepsilon),\label{eq:two-scale-g}\\
 \frac{D_m(V_m)}{N_m}&=10-2\cos\varepsilon+2\sin\varepsilon
            +\frac{4(1-\cos\varepsilon)}{\varepsilon^2}
             =10+O(\varepsilon),\label{eq:two-scale-D}\\
 \frac{T_m[V_m,V_m,V_m]}{N_m}&=-\frac32+O(\varepsilon).
 \label{eq:two-scale-T}
\end{align}
All implicit constants are universal.
\end{proposition}

\begin{proof}
The map $p\mapsto(\ell\cdot p\bmod4,h\cdot p\bmod m)$ is onto $\Z/4\Z\times\Z/m\Z$: integer points $(0,a,b)$ realize arbitrary pairs. Hence site averages are successive arithmetic means over the fast and slow phases. Every local increment of $\ell\cdot p$ and $h\cdot p$ lies in $\{-1,0,1\}$. Each cosine difference is at most two, and each normalized sine difference is at most one, proving \eqref{eq:two-scale-eta}.

The three distinct pairs of characters have frequencies $q\ell$, $q\ell+\varepsilon h$, $\varepsilon h$. The first and third are equality modes; for the second, substitution in the explicit symbol gives
\[
 c_1c_2c_3=\frac{1+\cos\varepsilon-\sin\varepsilon}{4},
 \qquad v^TN((q\ell+\varepsilon h)/2)v=\frac{57}{2}(1-\cos\varepsilon).
\]
Each real cosine or sine contributes half its complex quadratic value. These identities prove \eqref{eq:two-scale-g}, while summing the scalar bond symbol over the three modes proves \eqref{eq:two-scale-D}.

For the cubic coefficient, write $d_{p,r}=V_m(p+r)-V_m(p)$ and freeze the slow phase:
\[
 \widetilde d_{p,r}
 =v\{F_p(\ell\cdot r)-F_p(0)\}+w\cos\theta_p(h\cdot r),
 \quad F_p(s)=\cos(\psi_p+qs)+\cos(\psi_p+qs+\theta_p).
\]
The local arrays $d_p,\widetilde d_p$ are uniformly bounded and their difference is $O(\varepsilon)$ by scalar Taylor estimates. Thus the cubic polynomial $P(d)=D^3\Phi(0)[d,d,d]/6$ satisfies $P(d_p)=P(\widetilde d_p)+O(\varepsilon)$ uniformly. Lemma~\ref{lem:local-cubic} applies with
\[
 a_p=-\sin\psi_p-\sin(\psi_p+\theta_p),\quad
 b_p=-\cos\psi_p-\cos(\psi_p+\theta_p),\quad c_p=\cos\theta_p.
\]
Here $a_p^2+b_p^2=2+2\cos\theta_p$, while the fast mean of $a_p$ and slow mean of $\cos\theta_p$ are zero, and the slow mean of $\cos^2\theta_p$ is $1/2$. Equation~\eqref{eq:local-cubic} therefore gives $N_m^{-1}\sum_pP(\widetilde d_p)=-1/4$, proving \eqref{eq:two-scale-T}. Only derivatives at zero are used, so the unscaled field need not lie in the nonlinear neighborhood.
\end{proof}

For the remaining arguments, split
\[
 s_m=v\cos\psi+\frac{w}{\varepsilon}\sin\theta\in K_m,
 \qquad z_m=v\cos(\psi+\theta),\qquad V_m=s_m+z_m,
\]
and put $A=g_m(z_m)/N_m$, $B=T_m[s_m,s_m,z_m]/N_m$. The preceding proposition gives
\[
 A=\frac{19(1-\cos\varepsilon)}{128}\sim\frac{19}{256}\varepsilon^2,
 \qquad B=-\frac12+O(\varepsilon).
\]
Indeed, the pure $s_m$ cubic vanishes by Proposition~\ref{prop:cubic}, and all terms with at least two $z_m$ factors vanish by resonance. In fast and slow indices, $s_m$ has modes $(\pm1,0),(0,\pm1)$, and $z_m$ has modes $\pm(1,1)$. Three $z_m$ modes or one fast $s_m$ mode with two $z_m$ modes have odd fast sum; one slow $s_m$ mode with two $z_m$ modes leaves a nonzero slow index. Thus $T_m[V_m,V_m,V_m]=3N_mB$.

The same bond calculation gives
\[
 \frac{D_m(s_m)}{N_m}=4+\frac{4(1-\cos\varepsilon)}{\varepsilon^2}\leq6,
 \quad \frac{D_m(z_m)}{N_m}\leq16,\quad
 \eta(s_m)^2=2+2c_\varepsilon+2c_\varepsilon^2\leq6,
 \quad \eta(z_m)\leq2\sqrt2,
\]
where $c_\varepsilon=\sin\varepsilon/\varepsilon$. For the exact maximum in $\eta(s_m)$, the fast bond difference is $\pm1$, the largest normalized slow difference is $c_\varepsilon$, and independent phase indices attain matching signs on a bond with both increments one; use $v\cdot w=1$. In particular,
$d_m:=\eta(s_m)^2D_m(s_m)/N_m\to36$.

Taking $u_m=\varepsilon V_m$ in \eqref{eq:fourth-remainder} also gives
\begin{equation}\label{eq:linear-obstruction}
 \frac{E_m(u_m)-E_m(0)-D_m(u_m)/192}{N_m}
 =-\frac14\varepsilon^3+O(\varepsilon^4),\qquad
 \frac{E_m(u_m)-E_m(0)}{D_m(u_m)}
 =\frac1{192}-\frac{\varepsilon}{40}+O(\varepsilon^2).
\end{equation}

\begin{proof}[Proof of Theorem~\ref{thm:necessary-growth}]
For fixed $m$, take $u(t)=ts_m+\lambda t^2z_m$. Since $s_m\in\ker g_m$ and its pure cubic vanishes, the coefficient of $t^4$ in the normalized defect from $D_m/192$ is
\[
 \frac12A\lambda^2+\frac12B\lambda+F_m,\qquad
 F_m=\frac{D^4E_m(0)[s_m,s_m,s_m,s_m]}{24N_m},
 \qquad |F_m|\leq K_4d_m.
\]
Its first two terms are minimized at $\lambda=-B/(2A)$, with value $-B^2/(8A)$. For this fixed $m$, $\lambda$ is finite, so $t\to0$ enters every prescribed positive neighborhood. Moreover, the finite maximum defining $\eta$ is continuous, and
\[
 \frac{\eta(u(t))^2D_m(u(t))}{N_mt^4}\longrightarrow d_m.
\]
The assumed inequality implies $A_m\geq B^2/(8Ad_m)-K_4$. Multiplying by $\varepsilon^2$ and using $A/\varepsilon^2\to19/256$, $B^2\to1/4$, $d_m\to36$ proves \eqref{eq:necessary-growth}. Constants within $1$ of the infimum defining $C_m^{\rm loc}$ give the final assertion.
\end{proof}

\begin{proof}[Completion of the proof of Theorem~\ref{thm:crossover}]
The growth bounds follow from Theorem~\ref{thm:necessary-growth} and \eqref{eq:matching}. For the uniform lower loss, choose
\[
 u=t(s_m+a z_m),\qquad a=\min\{t/\varepsilon^2,1\}.
\]
The bounds above give $\eta(u)\leq6t$ and $D_m(u)\leq22N_mt^2$. For all sufficiently large $m$, $B\leq-1/4$ and $A\leq19\varepsilon^2/256$. The cubic cancellations reduce \eqref{eq:fourth-remainder} to
\[
 \frac{E_m(u)-E_m(0)-D_m(u)/192}{N_m}
 =\frac12Aa^2t^2+\frac12Bat^3+\widetilde R_4,
 \qquad |\widetilde R_4|\leq792K_4t^4.
\]
The two leading terms are at most $-45\min\{t^4/\varepsilon^2,t^3\}/512$. Choose $m_0$ sufficiently large and $t_0>0$ sufficiently small that
\[
 \varepsilon^2\leq\frac{45}{1024\cdot792K_4}\quad(m\geq m_0),
 \qquad t_0\leq\frac{45}{1024\cdot792K_4},
\]
and $6t_0$ lies in the uniform analytic radius. These choices absorb the remainder in both ranges $t\leq\varepsilon^2$ and $t\geq\varepsilon^2$. Set $t=\rho/6$ and $\rho_0=6t_0$. Dividing the resulting negative defect by $D_m(u)\leq22N_mt^2$ gives
\begin{equation}\label{eq:crossover}
 \frac1{192}-c_m(\rho)
 \geq\frac{45}{22528}\min\left\{\frac{m^2\rho^2}{144\pi^2},\frac\rho6\right\}
 \geq\frac{45}{3244032\pi^2}\min\{m^2\rho^2,\rho\}.
\end{equation}
Together with the upper estimate, this proves Theorem~\ref{thm:crossover}.
\end{proof}

\begin{proof}[Proof of the all-period assertion in Theorem~\ref{thm:loss}]
Equation~\eqref{eq:cubic-remainder} gives the upper loss bound. For the reverse bound, put $t=\rho/(5\sqrt2)$ and choose $4\mid m$ so large that $0<\varepsilon\leq t$. Then $tV_m$ is admissible, and \eqref{eq:two-scale-g}--\eqref{eq:two-scale-T}, with \eqref{eq:fourth-remainder}, give uniformly
\[
 \frac{E_m(tV_m)-E_m(0)}{D_m(tV_m)}
 =\frac1{192}-\frac{t}{40}+O(t^2).
\]
For sufficiently small $\rho>0$, the loss is at least $t/80$. Hence
\[
 \frac{\rho}{400\sqrt2}\leq\frac1{192}-c_{\rm all}(\rho)
 \leq\frac{10\sqrt2M_3}{3}\rho,
\]
as required.
\end{proof}

\section{Global optimality for two-site rectangular periods}\label{sec:global}
We now prove Theorem~\ref{thm:global} for two points per rectangular period cell, allowing arbitrary aspect ratios and relative shifts. The resulting point set need not be a Bravais lattice, and the proof permits changes of Voronoi incidence.

Let $\ell=(\ell_1,\ell_2,\ell_3)\in(0,\infty)^3$ satisfy $\ell_1\ell_2\ell_3=1$, and put
\[
 D_\ell=\diag(\ell_1,\ell_2,\ell_3),\qquad
 \Omega_\ell=D_\ell[0,1)^3.
\]
For $a\in[0,1]^3$, define
\begin{equation}\label{eq:global-energy}
 L_\ell(a)=D_\ell\bigl(\Z^3\cup(\Z^3+a)\bigr),\qquad
 Q(\ell,a)=\int_{\Omega_\ell}\dist(x,L_\ell(a))^2\,dx.
\end{equation}
Boundary values of the shift are identified modulo $\Z^3$. Coincident cosets can be included by continuity. The period cell has volume one, so the BCC value in this normalization is $19/128$.

\begin{proposition}[Global rectangular-period stability]\label{prop:global-rectangular}
Every configuration \eqref{eq:global-energy} satisfies
\[
 Q(\ell,a)\geq\frac{19}{128}.
\]
Equality holds exactly when $\ell_1=\ell_2=\ell_3=1$ and $a=(1/2,1/2,1/2)$ modulo $\Z^3$.
More precisely, relabel the axes so that $\ell_1=\max_i\ell_i$, and set
\begin{equation}\label{eq:global-shape-profile}
 f(t)=
 \begin{cases}
 t^{-2/3}\left(\dfrac1{48}+\dfrac t6-\dfrac{t^2}{24}\right),
       &0<t\leq\dfrac12,\\[2mm]
 t^{-2/3}\left(\dfrac{1+t}{12}-\dfrac1{48t}+\dfrac1{384t^2}\right),
       &\dfrac12\leq t\leq1.
 \end{cases}
\end{equation}
Then
\begin{equation}\label{eq:global-defect}
 Q(\ell,a)-\frac{19}{128}
 \geq f(\ell_1^{-3})-\frac{19}{128}
 +\frac{(\ell_2-\ell_3)^2}{96}
 +\frac{\displaystyle\sum_{i=1}^3\ell_i^4(a_i-\tfrac12)^2}
        {\displaystyle12\sum_{i=1}^3\ell_i^2}.
\end{equation}
The function $f$ is strictly decreasing on $(0,1]$, with $f(1)=19/128$ and $f(t)\sim t^{-2/3}/48$ as $t\downarrow0$.
\end{proposition}

The proof of Proposition~\ref{prop:global-rectangular} has two steps. An exact integral identity first shows that centering the second point lowers the energy for every fixed aspect ratio. A two-variable reduction then gives the minimum of the centered energy. Both steps are global and use no local expansion.

\begin{lemma}[Exact rectangular reduction]\label{lem:box-reduction}
For $w=(w_1,w_2,w_3)\in[0,\infty)^3$, let
\[
 A(w)=\frac18\int_{[-1,1]^3}\left|\sum_{i=1}^3w_iv_i\right|\,dv.
\]
With $w_i=\ell_i^2a_i(1-a_i)$, one has
\begin{equation}\label{eq:box-reduction}
 Q(\ell,a)=\frac1{12}\sum_{i=1}^3\ell_i^2-\frac12A(w).
\end{equation}
\end{lemma}
\begin{proof}
Write $d(s)=\dist(s,\Z)$. For $0<a<1$, the periodic function
$d(s)^2-d(s-a)^2$ has two affine branches, with slopes $2a$ and
$-2(1-a)$ and common range $[-a(1-a),a(1-a)]$. Changing variables
on these branches gives, for every continuous $g$,
\[
 \int_0^1g\bigl(d(s)^2-d(s-a)^2\bigr)\,ds
 =\frac12\int_{-1}^1g\bigl(a(1-a)v\bigr)\,dv.
\]
The identity also holds at $a=0,1$, when both sides equal $g(0)$.

For $x=D_\ell s$, the squared distances to the two cosets are
$\sum_i\ell_i^2d(s_i)^2$ and
$\sum_i\ell_i^2d(s_i-a_i)^2$. Integrating
$\min\{u,v\}=(u+v-|u-v|)/2$, using
$\int_0^1d(s)^2\,ds=1/12$, and applying the preceding identity
successively in the three coordinates yields
\[
 Q(\ell,a)=\frac1{12}\sum_i\ell_i^2
 -\frac1{16}\int_{[-1,1]^3}\left|\sum_iw_iv_i\right|\,dv.
\]
This is \eqref{eq:box-reduction}.
\end{proof}

\begin{lemma}[Global centering estimate]\label{lem:global-centering}
For every $\ell$ and $a$ as above,
\begin{equation}\label{eq:global-centering}
 Q(\ell,a)-Q\bigl(\ell,(\tfrac12,\tfrac12,\tfrac12)\bigr)
 \geq
 \frac{\displaystyle\sum_i\ell_i^4(a_i-\tfrac12)^2}
      {\displaystyle12\sum_i\ell_i^2}.
\end{equation}
\end{lemma}
\begin{proof}
We first establish the width estimate
\begin{equation}\label{eq:width-derivative}
 \frac{\partial A}{\partial w_i}(w)
 \geq\frac{w_i}{3\sum_jw_j},\qquad w_i>0.
\end{equation}
Let $b,c$ be the other widths. The normalized integral of
$g(bv+cz)$ over $[-1,1]^2$ has an even density that is
nonincreasing on $[0,\infty)$ and supported on $[-b-c,b+c]$:
for $b,c>0$, this density is
\[
 \frac{\bigl|[-b,b]\cap[y-c,y+c]\bigr|}{4bc}.
\]
Its average on a centered subinterval is therefore at least its
average on the full support. Consequently,
\[
 H(s):=\frac14\int_{\substack{(v,z)\in[-1,1]^2\\|bv+cz|\leq s}}\,dv\,dz
 \geq\min\left\{\frac{s}{b+c},1\right\},\qquad b+c>0.
\]
The same bound holds when one width vanishes, by the constant
density on the remaining interval. If both vanish, $H(s)=1$
for $s>0$. Thus in every case $H(s)\geq s/(w_i+b+c)$ for
$0\leq s\leq w_i$. Integrating first in the variable of width $w_i$ gives
\[
 A_i(w)=\frac18\int_{[-1,1]^2}
   \left(1-\frac{(bv+cz)^2}{w_i^2}\right)_+\,dv\,dz
 =\int_0^{w_i}\frac{s}{w_i^2}H(s)\,ds
 \geq\frac{w_i}{3\sum_jw_j},
\]
where $r_+=\max\{r,0\}$.

Put $W_i=\ell_i^2/4$ and integrate \eqref{eq:width-derivative}
along the segment from $w$ to $W$. Since every point of this
segment has width sum at most $\sum_iW_i$,
\[
 A(W)-A(w)\geq\frac{\sum_i(W_i^2-w_i^2)}{6\sum_iW_i}.
\]
For $\delta_i=a_i-1/2$, one has
\[
 W_i^2-w_i^2
 =\ell_i^4\delta_i^2\left(\frac12-\delta_i^2\right)
 \geq\frac14\ell_i^4\delta_i^2.
\]
Combining these inequalities with \eqref{eq:box-reduction} and
$\sum_iW_i=\frac14\sum_i\ell_i^2$ proves the claim.
Zero endpoints follow by continuity.
\end{proof}

\begin{proof}[Proof of Proposition~\ref{prop:global-rectangular}]
It remains to estimate the centered energy. For $x\geq\max\{y,z\}>0$, three successive integrations, using $s|s|^3/24$ as a third antiderivative of $|s|$, give
\begin{align*}
 A(x,y,z)
 &=\frac{(x+y+z)^4-(x+y-z)^4-(x-y+z)^4
       +(x-y-z)|x-y-z|^3}{96xyz}\\
 &=\frac x2+\frac{y^2+z^2}{6x}
       -\frac{(y+z-x)_+^4}{48xyz}.
\end{align*}
Take $x=\ell_1^2$, $y=\ell_2^2$, $z=\ell_3^2$. Since $xyz=1$ and $A$ is homogeneous of degree one, \eqref{eq:box-reduction} becomes
\begin{equation}\label{eq:centered-rectangular-formula}
 Q\bigl(\ell,(\tfrac12,\tfrac12,\tfrac12)\bigr)
 =\frac x{48}+\frac{y+z}{12}
  -\frac{y^2+z^2}{48x}
  +\frac{(y+z-x)_+^4}{384xyz}.
\end{equation}
This formula covers both possible signs of $y+z-x$.

Set $t=x^{-3/2}=\ell_1^{-3}$ and $s=(y+z)/x$.
Since $x=\max\{x,y,z\}$ and $xyz=1$, one has
$0<t\leq1$ and $2t\leq s\leq2$. Formula
\eqref{eq:centered-rectangular-formula} becomes
\[
 384Q\bigl(\ell,(\tfrac12,\tfrac12,\tfrac12)\bigr)
 =xG(s,t),\qquad
 G(s,t)=8+32s-8s^2+16t^2+\frac{(s-1)_+^4}{t^2}.
\]
On this range, $\partial_sG\geq4$. For $s\leq1$ this follows
from $32-16s\geq16$; for $1\leq s\leq2$, use $t\leq1$ to get
\[
 \partial_sG(s,t)\geq32-16s+4(s-1)^3
 =4+4(2-s)\bigl(3-(s-1)-(s-1)^2\bigr)\geq4.
\]
Thus
\[
 Q\bigl(\ell,(\tfrac12,\tfrac12,\tfrac12)\bigr)
 \geq\frac{xG(2t,t)}{384}+\frac{x(s-2t)}{96}
 =f(t)+\frac{(\ell_2-\ell_3)^2}{96}.
\]
Here $x(s-2t)=(\sqrt y-\sqrt z)^2$, and evaluation of
$G(2t,t)$ gives the two expressions in
\eqref{eq:global-shape-profile}. Lemma~\ref{lem:global-centering}
now proves \eqref{eq:global-defect}.

The two expressions in \eqref{eq:global-shape-profile} agree at $t=1/2$, and direct differentiation gives
\[
 f'(t)=
 \begin{cases}
 -\dfrac{(2t-1)^2}{72t^{5/3}},&0<t\leq\dfrac12,\\[2mm]
 \dfrac{(t-1)(2t-1)^2}{144t^{11/3}},
      &\dfrac12\leq t\leq1.
 \end{cases}
\]
Thus $f$ is strictly decreasing, and $f(t)\geq f(1)=19/128$, with equality only at $t=1$. The stated asymptotic follows from its first expression. Finally, equality with the BCC energy forces $\ell_1=1$ and hence $\ell_2=\ell_3=1$, since $\ell_1$ is maximal and the product is one. The last term in \eqref{eq:global-defect} then forces $a=(1/2,1/2,1/2)$. The converse is the reference BCC identity.
\end{proof}

For the cubic period, the same representation gives a sharper global
estimate. We write $\mathbf 1=(1,1,1)$ and regard the shift as
an element of $\mathbb T^3=\mathbb R^3/\mathbb Z^3$.

\begin{proposition}[A sharp global estimate on the cubic torus]
\label{prop:global-cubic}
For every $a\in\mathbb T^3$, put
$\delta_i=\operatorname{dist}(a_i,\frac12+\mathbb Z)$.  Then
\begin{equation}\label{eq:global-cubic-gap}
 Q(\mathbf1,a)-\frac{19}{128}
 \geq \frac7{96}\sum_{i=1}^3\delta_i^2
       +\frac1{24}\sum_{1\leq i<j\leq3}\delta_i^2\delta_j^2
       +\frac52\delta_1^2\delta_2^2\delta_3^2.
\end{equation}
In particular, $7/96$ is the largest constant $c$ for which
$Q(\mathbf1,a)-19/128\geq c\sum_i\delta_i^2$ holds for every shift.
\end{proposition}

\begin{proof}
Set $x_i=1-4\delta_i^2\in[0,1]$. By
\eqref{eq:box-reduction},
$Q(\mathbf1,a)=1/4-A(x/4)/2$.
The function $A$ is convex in each coordinate, and its vertex values are
\[
 A(\varepsilon/4)=\frac14b_{\varepsilon_1+\varepsilon_2+\varepsilon_3},
 \qquad \varepsilon\in\{0,1\}^3,\qquad
 (b_0,b_1,b_2,b_3)=\left(0,\frac12,\frac23,\frac{13}{16}\right).
\]
Applying the affine chord bound successively in the three coordinates gives
\[
 A(x/4)\leq\frac14\left(
 \frac12\sum_i x_i-\frac13\sum_{i<j}x_ix_j
 +\frac5{16}x_1x_2x_3\right).
\]
Substituting $x_i=1-4\delta_i^2$ proves
\eqref{eq:global-cubic-gap}. At $a=(0,1/2,1/2)$,
$Q(\mathbf1,a)=1/6$ and $\sum_i\delta_i^2=1/4$, while
the last two terms vanish. Hence $7/96$ is attained and is optimal.
\end{proof}

We conclude this section by classifying all critical shifts on the
cubic torus. The classification follows from the regularity and
monotonicity properties of the same integral representation.

\begin{proposition}[All critical shifts on the cubic torus]
\label{prop:global-critical}
The function $a\mapsto Q(\mathbf1,a)$ belongs to
$C^2(\mathbb T^3\setminus\{0\})$.  Its critical points are precisely
the seven nonzero elements of $\{0,\frac12\}^3$, modulo
$\mathbb Z^3$.  The BCC shift $(\frac12,\frac12,\frac12)$ is the
unique global minimum.  Up to a coordinate permutation, the other
critical points and their Hessians are
\begin{align*}
 D_a^2Q\left(\mathbf1,(\tfrac12,\tfrac12,\tfrac12)\right)
       &=\frac{13}{48}I_3,\\
 D_a^2Q\left(\mathbf1,(\tfrac12,\tfrac12,0)\right)
       &=\operatorname{diag}\left(\frac13,\frac13,-\frac23\right),\\
 D_a^2Q\left(\mathbf1,(\tfrac12,0,0)\right)
       &=\operatorname{diag}\left(\frac12,-\frac23,-\frac23\right).
\end{align*}
Consequently the three shifts of the second type have Morse index
one, and the three shifts of the third type have Morse index two.
The coincident shift $a=0$ is the unique global maximum, with value
$1/4$, and the energy is not differentiable there.
\end{proposition}

\begin{proof}
Extend $A$ to $\mathbb R^3$ by its defining integral. It is even in
each coordinate and belongs to $C^2(\mathbb R^3\setminus\{0\})$.
Indeed, near a point with $w_1>0$, integrate first in the first variable:
\[
 \begin{aligned}
 B(b,s)&=\frac12\int_{-1}^1|bt+s|\,dt
 =\begin{cases}
 |s|,&|s|\geq b,\\[2pt]
 (b^2+s^2)/(2b),&|s|\leq b,
 \end{cases}\\
 A(w)&=\frac14\int_{[-1,1]^2}B(w_1,w_2v+w_3z)\,dv\,dz.
 \end{aligned}
\]
The kernel is $C^1$ with locally bounded piecewise second derivatives.
Its junction sets have area zero when $(w_2,w_3)\ne(0,0)$;
when $(w_2,w_3)=(0,0)$, the junctions are absent nearby.
Differentiation under the integral and dominated convergence prove
the assertion. Near an integer shift, $w_i=|a_i|-a_i^2$.
Since evenness gives $A_i=0$ and $A_{ij}=0$ for $j\ne i$ on
$w_i=0$, the chain rule shows that both first and second derivatives
of $A(w(a))$ agree across $a_i=0$. This proves the claimed
regularity on the punctured torus.

If $a_i\notin\{0,1/2\}$ modulo $\mathbb Z$, then $w_i>0$ and
\eqref{eq:width-derivative} gives
\[
 \partial_{a_i}Q=-\frac12(1-2a_i)A_i(w)\ne0.
\]
Conversely, a nonzero shift in $\{0,1/2\}^3$ is fixed modulo
$\mathbb Z^3$ by reflection in each coordinate, so all its first
derivatives and mixed second derivatives vanish.

Put $b=1/4$. At a half integer, the diagonal Hessian entry is
$A_i(w)$. Symmetry and degree-one homogeneity, applied to the
vertex values above, give $A_i=1/2$, $1/3$, and $13/48$ when
one, two, and three widths are active. At an integer, that entry is
$-A_{ii}(w)/2$. To compute $A_{ii}$ for an inactive width, integrate
first in an active variable using $B$. The remaining active width is
either $0$ or $b$, so $B_{ss}=1/b$ almost everywhere in the
remaining integral. Therefore
\[
 A_{ii}(w)=\frac1{4b}\int_{[-1,1]^2}v^2\,dv\,dz
 =\frac1{3b}=\frac43.
\]
These values prove the three Hessian formulas and their Morse indices.

The unique minimum follows from \eqref{eq:global-cubic-gap}.
Since $A(w)>0$ for $w\ne0$, the coincident shift is the unique
maximum, of value $1/4$. Finally,
\[
 Q(\mathbf1,(t,0,0))=\frac14-\frac{|t|}{4}+\frac{t^2}{4},
 \qquad |t|\leq\frac12,
\]
proves nondifferentiability there.
\end{proof}

\appendix
\section{Exact arithmetic lemmas for the determinant}\label{app:certificates}

We record the two exact coefficient identities used in Lemma~\ref{lem:matrix}. The formulas below specify every integer coefficient and its zero set.

For a concise expansion of \eqref{eq:determinant}, put
\[
 T=a^2+b^2+c^2,\qquad
 U=a^2b^2+a^2c^2+b^2c^2,\qquad V=a^2b^2c^2.
\]
Direct expansion of that identity gives
\begin{equation}\label{eq:compact-determinant}
\begin{split}
 F_\sigma(a,b,c)={}&64000-102440T+38464T^2-24T^3\\
 &+75430U-38038TU-52390V+29854TV\\
 &+\sigma abc\bigl(51442-75904T+846T^2+66076U-52900V\bigr).
\end{split}
\end{equation}
Write $f^\sigma_{rsu}=[a^rb^sc^u]F_\sigma$, where $[a^rb^sc^u]$ denotes coefficient extraction. Formula~\eqref{eq:compact-determinant} specifies these integers completely; exactly $34$ are nonzero, each exponent is at most six, and $r+s+u\le9$. We use the convention
\[
 \binom nj=0\quad\hbox{if }n<0\hbox{ or }j<0\hbox{ or }j>n.
\]

\begin{lemma}[Positive-branch coefficient certificate]\label{lem:positive-certificate-data}
For $i,j,q,m\ge0$ with $i+j+q+m=8$, define
\begin{equation}\label{eq:positive-coefficient-formula}
 b_{ijqm}=\sum_{r,s,u=0}^6 f^+_{rsu}
       \binom{9-r-s-u}{m+1}
       \binom{8-r-s-m}{q}
       \binom{8-r-m-q}{j}.
\end{equation}
These integers vanish exactly at
\[
 (i,j,q,m)\in
 \{(8,0,0,0),(7,1,0,0),(7,0,1,0),(7,0,0,1)\}.
\]
At each of the remaining $161$ indices,
\begin{equation}\label{eq:positive-coefficient-bounds}
 16896\le b_{ijqm}\le123151676.
\end{equation}
In particular, $b_{0008}=64000$. The polynomial
\[
 B(t,h,k,l)=\sum_{i+j+q+m=8}b_{ijqm}t^ih^jk^ql^m
\]
satisfies the exact identity \eqref{eq:positive-certificate}.
\end{lemma}

\begin{proof}
Set $L=t+h+k+l$. The contribution of the monomial $a^rb^sc^u$ to
$L^9F_+(t/L,(t+h)/L,(t+h+k)/L)$ is
\[
 f^+_{rsu}\,t^r(t+h)^s(t+h+k)^uL^{9-r-s-u}.
\]
Selecting the power $l^{m+1}$ first leaves
\[
 f^+_{rsu}\binom{9-r-s-u}{m+1}
 t^r(t+h)^s(t+h+k)^{8-r-s-m}.
\]
Selecting $k^q$ and then $h^j$ gives precisely
\eqref{eq:positive-coefficient-formula}; the zero-binomial convention includes every out-of-range case. The term with $l=0$ factors as
\[
 6h^2(t+h+k)^3
 (17h^2+32hk+16ht+16k^2+32kt)
 (119h^2+246hk+299ht+123k^2+246kt+176t^2),
\]
which proves the asserted polynomial identity.

For the bounds, substitute the $34$ nonzero coefficients of
\eqref{eq:compact-determinant} into \eqref{eq:positive-coefficient-formula}.
The $\binom{11}{3}=165$ indices $i+j+q+m=8$ give exactly the four stated zeros, minimum positive value $b_{6200}=16896$, and maximum $b_{2222}=123151676$. This proves \eqref{eq:positive-coefficient-bounds} by finite integer arithmetic.
\end{proof}

\begin{lemma}[Negative-branch coefficient certificate]\label{lem:negative-certificate-data}
For $0\le i,j,k\le6$, define
\begin{equation}\label{eq:negative-coefficient-formula}
 Z^-_{ijk}=\sum_{r,s,u=0}^6 f^-_{rsu}
       \binom{6-r}{i-r}\binom{6-s}{j-s}\binom{6-u}{k-u}.
\end{equation}
These integers vanish exactly at the coordinate permutations of
$(0,0,6)$ and $(0,1,6)$. At the remaining $334$ positions,
\begin{equation}\label{eq:negative-coefficient-bounds}
 11808\le Z^-_{ijk}\le147112212.
\end{equation}
In particular, $Z^-_{000}=64000$, and these coefficients satisfy the exact identity \eqref{eq:negative-certificate}.
\end{lemma}

\begin{proof}
For $0\le r\le6$, the binomial theorem gives
\[
 a^r=\sum_{i=r}^6\binom{6-r}{i-r}a^i(1-a)^{6-i}.
\]
Apply this identity separately in the three variables to
\eqref{eq:compact-determinant} with $\sigma=-1$. The result is
\eqref{eq:negative-certificate}, with coefficients given by
\eqref{eq:negative-coefficient-formula}. No binomial factors are implicit in that basis.

For the bounds, symmetry reduces the calculation to the $84$ indices
$0\le i\le j\le k\le6$. Substitution of the $34$ source coefficients from
\eqref{eq:compact-determinant} gives exactly the stated zeros, minimum positive value $Z^-_{026}=11808$, and maximum $Z^-_{222}=147112212$. Coordinate permutations give all $343$ coefficients and prove \eqref{eq:negative-coefficient-bounds}.
\end{proof}

Thus both determinant inequalities hold on the entire closed parameter cube by exact integer arithmetic.

\section{Exact local cell identities}\label{app:cells}
We give a finite rational construction for the cell moment used in Section~\ref{sec:fixed}, and for the cubic coefficient used in the growing-period argument. All vertices, face triangulations, and moment formulas are specified below, so the identities can be checked directly by rational algebra.

For $i\ne j$, $\sigma,\tau\in\{-1,1\}$, write
\[
 q_{ij}^{\sigma\tau}(0)=\frac\sigma2e_i+\frac\tau4e_j,
\]
and let $k$ be the remaining coordinate. These are the $24$ reference vertices. Their three active normals are
\begin{equation}\label{eq:active-triples}
 r_1=\sigma e_i,\qquad
 r_2=\tfrac12(\sigma e_i+\tau e_j+e_k),\qquad
 r_3=\tfrac12(\sigma e_i+\tau e_j-e_k).
\end{equation}
For general local differences $d=(d_r)_{r\in\mathcal N}$, put $n_r=r+d_r$ and continue the vertex by
\begin{equation}\label{eq:rational-vertex}
 q_{ij}^{\sigma\tau}(d)
 =A^{-1}\begin{pmatrix}|n_{r_1}|^2/2\\|n_{r_2}|^2/2\\|n_{r_3}|^2/2\end{pmatrix},
 \qquad A=\begin{pmatrix}n_{r_1}^T\\n_{r_2}^T\\n_{r_3}^T\end{pmatrix}.
\end{equation}
The six axis faces have the cyclic vertex lists
\[
 \bigl(q_{ij}^{\sigma,+},q_{ik}^{\sigma,+},
 q_{ij}^{\sigma,-},q_{ik}^{\sigma,-}\bigr),
 \qquad \{j,k\}=\{1,2,3\}\setminus\{i\},\quad j<k.
\]
The face with normal $(\sigma_1,\sigma_2,\sigma_3)/2$ has the cyclic list
\[
 \bigl(q_{12}^{\sigma_1\sigma_2},q_{21}^{\sigma_2\sigma_1},
 q_{23}^{\sigma_2\sigma_3},q_{32}^{\sigma_3\sigma_2},
 q_{31}^{\sigma_3\sigma_1},q_{13}^{\sigma_1\sigma_3}\bigr).
\]
On each face with ordered list $(q_1,\ldots,q_s)$, take the triangles
$(q_1,q_j,q_{j+1})$, $2\leq j\leq s-1$, and cone them to the origin.
There are $6\cdot2+8\cdot4=44$ tetrahedra. For each ordered triple, let $\epsilon$ be the sign of its reference determinant. In the neighborhood of Lemma~\ref{lem:local-cell} the sign is unchanged, and the tetrahedron $T=(0,q_1,q_2,q_3)$ has
\begin{equation}\label{eq:tetrahedron-moment}
 |T|=\frac{\epsilon\det(q_1,q_2,q_3)}6,
 \qquad
 \int_T|y|^2\,dy
 =\frac{\epsilon\det(q_1,q_2,q_3)}{120}
 \left(\sum_{i=1}^3|q_i|^2+|q_1+q_2+q_3|^2\right).
\end{equation}
The formula follows by integrating the quadratic monomials in barycentric coordinates. Summing these $44$ expressions gives the full volume and cell moment.

For the phase $a=0$ of~\eqref{eq:four-phase-field}, introduce the independent parameters $t,b$ and set $v=(1,-1,0)$ and $h_0=(1,1,0)$. Its deformed neighbor normals are
\begin{equation}\label{eq:fourphase-normals}
 n_r=\begin{cases}
 r-2b h_0,&r_1+r_2=1,\\
 r-2tv-2b h_0,&r_1+r_2=-1,\\
 r,&r_1+r_2=0.
 \end{cases}
\end{equation}
Here $r_1,r_2$ denote the first two coordinates of $r$.
Substitution of~\eqref{eq:fourphase-normals} into~\eqref{eq:rational-vertex} and~\eqref{eq:tetrahedron-moment}, followed by collection over the common denominator, gives the two exact identities
\begin{equation}\label{eq:fourphase-rational-sums}
 |K|=\frac12,\qquad
 \int_K|y|^2\,dy
 =\frac{57+32t^2-192b^2-128bt^2-12288b^4}{768(1-16b^2)}.
\end{equation}
The corresponding determinant sum is
$\sum\epsilon\det(q_1,q_2,q_3)=3$; the moment numerator is obtained from
the same $44$ explicitly specified summands.
Putting $b=\lambda t^2$ gives~\eqref{eq:fourphase-energy}.
All four phase cells have the same volume and moment: the coordinate maps
\[
 P(x_1,x_2,x_3)=(x_2,x_1,x_3),\qquad
 S(x_1,x_2,x_3)=(-x_2,-x_1,x_3)
\]
satisfy $u_{Pp+e_1+e_2}=Pu_p$ and $u_{Sp+e_1}=Su_p$.
The induced isometries preserve the complete point set and act transitively on the four phases.

\begin{lemma}[Certified geometry of the four-phase cells]\label{lem:fourphase-geometry}
Set $b=4t^2/41$ in~\eqref{eq:fourphase-normals}, and let $q$ range over
the $24$ vertices defined by~\eqref{eq:rational-vertex}. For
$|t|\leq1/1000$, every active normal matrix $A$ and every inactive
neighbor normal $n_r$ at $q$ satisfy
\[
 |\det A|\geq\frac{49}{100},\qquad
 \frac{|n_r|^2}{2}-n_r\cdot q\geq\frac6{25}.
\]
For every tetrahedron in the specified face triangulations,
\[
 \epsilon\det(q_1,q_2,q_3)\geq\frac1{25}.
\]
In addition, every axis normal and every listed vertex satisfy
\[
 \|n_{\pm e_i}-(\pm e_i)\|_1<\frac1{200},\qquad
 |q|^2<\frac8{25}.
\]
Consequently these vertices and the $44$ tetrahedra describe the entire
fourteen-plane intersection throughout this interval.
\end{lemma}
\begin{proof}
We compare each vertex with its reference position $q_0$.
Every normal perturbation $d_r=n_r-r$ has at most two nonzero coordinates,
and~\eqref{eq:fourphase-normals} gives
\[
 \|d_r\|_\infty\leq\varepsilon:=\frac{2001}{10^6}.
\]
For every active triple, $\|A_0^{-1}\|_\infty=3$ and
$\|A-A_0\|_\infty\leq2\varepsilon$. Since
$|r_i-q_0|_1=3/4$, the vertex equations give
\[
 \|q-q_0\|_\infty
 \leq\frac{3(3\varepsilon/4+\varepsilon^2)}{1-6\varepsilon}
 <\delta:=\frac{23}{5000}.
\]
Here the numerator bounds
$Aq-Aq_0=((r_i-q_0)\cdot d_{r_i}+|d_{r_i}|^2/2)_{i=1}^3$.
Also, direct substitution in the active triples gives
\[
 |\det A|\geq\frac12-\frac8{41}t^2>\frac{49}{100}.
\]
At an inactive incidence, put
$s_0=|r|^2/2-r\cdot q_0$. The reference vertices give $s_0\geq1/4$;
if $s_0=1/4$, then $|r-q_0|_1=5/4$, and otherwise $s_0\geq1/2$
and $|r-q_0|_1\leq9/4$. As $|r|_1\leq3/2$, expanding the perturbed
slack yields, in both cases,
\[
 \frac{|n_r|^2}{2}-n_r\cdot q
 \geq\min\left\{
 \frac14-\frac32\delta-\frac54\varepsilon-3\varepsilon\delta,
 \frac12-\frac32\delta-\frac94\varepsilon-3\varepsilon\delta
 \right\}>\frac6{25}.
\]
For every reference triangle, the displayed face lists give
\[
 |\det(q_1^0,q_2^0,q_3^0)|\geq\frac3{64},\qquad
 \sum_{\rm cyc}|q_i^0\times q_j^0|_1\leq\frac54,
 \qquad |q_i^0|_1=\frac34.
\]
By multilinearity of the determinant,
\[
 \epsilon\det(q_1,q_2,q_3)
 \geq\frac3{64}-\frac54\delta-\frac92\delta^2-6\delta^3
 >\frac1{25}.
\]
The quadratic error term follows from
$|\det(a,b,q_i^0)|\leq2\delta^2|q_i^0|_1$ when
$\|a\|_\infty,\|b\|_\infty\leq\delta$; the cubic term follows by expanding the six determinant monomials.
Finally,
\[
 |q|^2\leq\frac5{16}+\frac32\delta+3\delta^2<\frac8{25},
 \qquad
 \|n_{\pm e_i}-(\pm e_i)\|_1
 \leq4\cdot10^{-3}+\frac{16}{41}\,10^{-6}
 =\frac{10251}{2562500}<\frac1{200}.
\]
The axis estimates make the intersection bounded. Every listed vertex
is simple, and the strict inactive slacks preserve all its incident
edges. Thus the listed vertices form a set closed under adjacency;
connectedness of the vertex graph excludes further vertices.
The positive oriented determinants give the asserted triangulation.
\end{proof}

To exclude omitted generators, note that every vector in
$\Lambda\setminus(\mathcal N\cup\{0\})$ has length at least $\sqrt2$. Since
$|u_p(t)|\leq\sqrt2(|t|+\lambda t^2)$, each omitted deformed generator has squared length at least
\[
 2(1-2\cdot10^{-3}-2\lambda\cdot10^{-6})^2
 =\frac{26160657333001}{13132812500000}
 >\frac{19}{10}>4\frac8{25}.
\]
It lies beyond twice the covering radius, and its inequality is strictly redundant. Thus these are the actual Voronoi cells on the stated interval.
At $\lambda=4/41$, the companion identity
\begin{equation}\label{eq:fourphase-positive}
 \frac{E_\Gamma(u(t))-E_\Gamma(0)-D_\Gamma(u(t))/384}{N}
 =\frac{t^2(2825761-1048944t^2+430336t^4-188416t^6)}
 {80688(1681-256t^4)}
\end{equation}
is positive for $0<|t|\leq1/1000$, while~\eqref{eq:fourphase-defect} is negative. Indeed, the relevant numerator lower bounds are
$2825761-1048944\cdot10^{-6}-188416\cdot10^{-18}>0$ and
$68921-53792\cdot10^{-6}>0$.
Hence the explicit family obeys
\[
 \frac1{384}<
 \frac{E_\Gamma(u(t))-E_\Gamma(0)}{D_\Gamma(u(t))}
 <\frac1{192}
 \qquad(0<|t|\leq1/1000).
\]

\begin{lemma}[A local cubic identity]\label{lem:local-cubic}
Put $\ell=(1,1,0)$, $v=(1,-1,0)$, $h=(1,0,1)$ and $w=(1,0,-1)$. For $r\in\mathcal N$, define
\[
 d_r=v\{a(\ell\cdot r)+b(\ell\cdot r)^2\}+wc(h\cdot r).
\]
Then the cubic Taylor coefficient of the local cell moment satisfies
\begin{equation}\label{eq:local-cubic}
 \frac16D^3\Phi(0)[d,d,d]
 =-\frac14c(a^2+b^2)-\frac14ac^2.
\end{equation}
\end{lemma}
\begin{proof}
At each active triple in~\eqref{eq:active-triples}, set
\[
 A_0=(r_i^T)_{i=1}^3,\qquad A_1=(d_{r_i}^T)_{i=1}^3,
\]
and let the $i$th coordinates of $b_0,b_1,b_2$ be
$|r_i|^2/2$, $r_i\cdot d_{r_i}$ and $|d_{r_i}|^2/2$, respectively.
The coefficients through order three of the corresponding vertex of $K(sd)$ are
\begin{equation}\label{eq:vertex-cubic-recursion}
 \begin{aligned}
 x_0&=A_0^{-1}b_0,&
 x_1&=A_0^{-1}(b_1-A_1x_0),\\
 x_2&=A_0^{-1}(b_2-A_1x_1),&
 x_3&=-A_0^{-1}A_1x_2.
 \end{aligned}
\end{equation}
Insert $q(s)=x_0+sx_1+s^2x_2+s^3x_3$ into the $44$ tetrahedral
moment expressions~\eqref{eq:tetrahedron-moment}, and denote their summed
coefficient of $s^3$ by $P(a,b,c)$. The recurrence
\eqref{eq:vertex-cubic-recursion} gives, by coefficient collection,
\begin{equation}\label{eq:cubic-coefficient-check}
 \begin{aligned}
 P(a,b,0)&=0,\\
 P(a,0,c)&=-\tfrac14ac(a+c),\\
 P(0,b,c)&=-\tfrac14b^2c,
 \qquad [abc]P=0.
 \end{aligned}
\end{equation}
Here $[abc]P$ denotes the coefficient of $abc$. These exact polynomial
identities determine all ten coefficients of the homogeneous cubic:
the three plane restrictions determine every monomial except $abc$,
whose coefficient is given separately. Hence
$P=-\tfrac14a^2c-\tfrac14b^2c-\tfrac14ac^2$, proving
\eqref{eq:local-cubic}. All coefficient extractions in
\eqref{eq:cubic-coefficient-check} use the same explicit vertex recurrence
and tetrahedron sum, without a separate numerical input.
\end{proof}

\end{document}